\documentclass{amsart}

\usepackage{amssymb}

\let\cal\mathcal
\allowdisplaybreaks

\newtheorem{theorem}{Theorem}[section]
\newtheorem{lemma}[theorem]{Lemma}
\newtheorem{proposition}[theorem]{Proposition}
\newtheorem{corollary}[theorem]{Corollary}
\newtheorem{remark}[theorem]{Remark}

\theoremstyle{definition}
\newtheorem{definition}[theorem]{Definition}

\newcommand{\R}{\mathbb{R}}
\newcommand{\C}{\mathbb{C}}

\def\M{\cal{M}}
\def\T{\tau}
\def\E{\mathsf{E}}

\newcommand{\dem}{\noindent {\bf Proof. }}
\newcommand{\fin}{\hspace*{\fill} $\square$ \vskip0.2cm}

\begin{document}

\null

\

\

\

\title[Non-commutative Gundy's decomposition]
{Gundy's decomposition  for non-commutative martingales and
applications}

\author[Parcet and Randrianantoanina]
{Javier Parcet$^*$ and Narcisse Randrianantoanina$^\dag$}

\begin{abstract} We provide an analogue of Gundy's decomposition
for $L_1$-bounded  non-commutative martingales. An important
difference from the classical case is that for any $L_1$-bounded
non-commutative martingale, the decomposition consists of four
martingales. This is strongly related with the row/column nature
of non-commutative Hardy spaces of martingales. As applications,
we obtain simpler proofs of the weak type $(1,1)$ boundedness for
non-commutative martingale transforms and the non-commutative
analogue of Burkholder's weak type inequality for square
functions. A sequence $(x_n)_{n \ge 1}$ in a normed space
$\mathrm{X}$ is called 2-co-lacunary if there exists a bounded
linear map from the closed linear span of $(x_n)_{n \ge 1}$ to
$l_2$ taking each $x_n$ to the $n$-th vector basis of $l_2$. We
prove (using our decomposition) that any relatively weakly compact
martingale difference sequence in $L_1(\M,\T)$ whose sequence of
norms is bounded away from zero is 2-co-lacunary, generalizing a
result of Aldous and Fremlin to non-commutative $L_1$-spaces.
\end{abstract}

\address{Department of Mathematics, Universidad Aut\'{o}noma de Madrid,
28049, Spain}

\email{javier.parcet@uam.es}

\address{Department of Mathematics, Miami University, Oxford, OH
45056}

\email{randrin@muohio.edu}

\footnote{$*$ Partially supported by MTM2004-00678.}
\footnote{$\dag$ Partially supported by NSF DMS-0096696.}

\maketitle

\section*{Introduction}

The main motivation for this paper comes from a fundamental
decomposition of martingales due to Gundy \cite{G} which is
generally referred  to as the {\it Gundy's decomposition theorem}.
Gundy's theorem has been very useful in establishing weak type
$(1,1)$ boundedness of certain quasi-linear mappings such as
square functions and Doob's maximal functions. In particular,
certain classical inequalities such as the weak type $(1,1)$
boundedness of martingale transforms and Burkholder's weak type
inequality for square functions can be deduced from Gundy's
theorem. We refer to \cite{BG,G2,MT} for some variations of
Gundy's result and more applications and to Garsia's notes
\cite{GA} for a complete discussion on classical martingale
inequalities.

\vskip5pt

Gundy's decomposition theorem played a central role in classical
martingale theory and it can be regarded as a probabilistic
counterpart of the well known Calder\'{o}n-Zygmund decomposition for
integrable functions \cite{CZ} in harmonic analysis. Due to its
relevance in the classical theory, it is natural to consider
whether or not such decomposition theorem can be generalized to
the non-commutative setting. In this paper, we investigate
possible analogues of Gundy's theorem for non-commutative
martingales. We first recall this classical result:

\vskip5pt

\noindent \textbf{Gundy's decomposition theorem {\rm{\cite{G}}}.}
\label{Gundy-C} \emph{Let $f=(f_n)_{n \ge 1}$ be a martingale on a
probability space $(\Omega, \mathcal{F},\mathbb{P})$ that is
bounded in $L_1(\Omega,\mathbb{P})$, and $\lambda$  be a positive
real number.  Then there are three martingales $a$, $b$, and $c$
relative to the same filtration and satisfying the following
properties for some absolute constant $\mathrm{c}$:}
\begin{itemize}
\item[(i)] $f=a+b+c$;
\item[(ii)] the martingale $a$ satisfies $$\|a\|_1 \le \mathrm{c}
\|f\|_1, \quad \|a\|_2^2 \le \mathrm{c} \lambda\|f\|_1, \quad
\|a\|_\infty \le \mathrm{c} \lambda;$$
\item[(iii)] the martingale $b$ satisfies $$\sum_{k=1}^{\infty}
\|db_k\|_1 \le \mathrm{c} \|f\|_1;$$
\item[(iv)] the martingale $c$ satisfies
$$\lambda \mathbb{P} \Big( \Big\{ \sup_{k\geq 1}|dc_k|>0 \Big\}
\Big) \le \mathrm{c} \left\|f\right\|_1.$$
\end{itemize}

As a prominent subfield of the theory of non-commutative
probability,  the theory of non-commutative martingale
inequalities has achieved considerable progress in recent years.
Indeed, many classical inequalities have been reformulated  to
include non-commutative martingales. This general theme started
from  the fundamental paper of  Pisier and Xu \cite{PX1} where
they introduced non-commutative martingale Hardy spaces and
formulated the right analogue of Burkholder-Gundy inequalities. It
was their general functional analytic approach that led to the
renewed interests in this topic. Shortly after \cite{PX1}, Junge
obtained in \cite{J1} a non-commutative analogue of Doob's maximal
functions.  Extensions of Burkholder/Rosenthal inequalities for
conditioned square functions were proved by Junge and Xu in
\cite{JX}. Martingale BMO spaces were studied in \cite{PX1, M, JM}
and some weak type inequalities can be found in \cite{R,R3}. We
also refer the reader to  a recent survey by Xu \cite{X} for a
rather complete exposition of the subject.

\vskip5pt

Following this general theme, we analyze analogues of Gundy's
decomposition for non-commutative martingales. For this, we note
first that since the notion of supremum does not necessarily make
sense for a family of operators, we require an appropriate
reformulation of condition (iv) above. It is clear that the
following equalities of measurable sets hold: $$\Big\{ \sup_{k\geq
1} |dc_k| > 0 \Big\} = \bigcup_{k\geq 1} \Big\{ |dc_k| > 0 \Big\}
= \bigcup_{k\geq 1} \text{supp} |dc_k|.$$ That is, condition (iv)
is equivalent to: $$\lambda \mathbb{P} \Big( \bigcup_{k\geq 1}
\text{supp} |dc_k| \Big)\leq \mathrm{c} \|f\|_1.$$ A
non-commutative analogue of this condition can be formulated using
the notion of support projection of a measurable operator. Our
main result in this paper reads as follows:

\vskip5pt

\noindent \textbf{Theorem A.} \emph{Let $\cal{M}$ be a semifinite
von Neumann algebra equipped with a normal semifinite trace $\T$
and let $(\M_n)_{n \ge1}$ be an increasing filtration of von
Neumann subalgebras of $\M$. If $x=(x_n)_{n \ge 1}$ is a
$L_1$-bounded non-commutative martingale  with respect to
$(\M_n)_{n \ge 1}$ and $\lambda$ is a positive real number, there
exist four martingales $\alpha$, $\beta$, $\gamma$ and $\upsilon$
relative to the same filtration and satisfying:}
\begin{itemize}
\item[(i)] $x=\alpha +\beta + \gamma + \upsilon$;
\item[(ii)] the martingale $\alpha$ satisfies $$\|\alpha\|_1 \leq
\mathrm{c} \|x\|_1, \quad \|\alpha\|_2^2 \leq \mathrm{c}
\lambda\|x\|_1, \quad \|\alpha\|_\infty \leq \mathrm{c} \lambda;$$
\item[(iii)] the martingale $\beta$ satisfies $$\sum_{k=1}^{\infty}
\|d\beta_k\|_1 \le \mathrm{c} \|x\|_1;$$
\item[(iv)] $\gamma$ and $\upsilon$ are $L_1$-martingales with
$$\max \Big\{ \lambda \tau \Big( \bigvee_{k \ge 1} \mathrm{supp}
|d\gamma_k| \Big), \, \lambda \tau \Big( \bigvee_{k \ge 1}
\mathrm{supp} |d\upsilon_k^*| \Big) \Big\} \le \mathrm{c}
\|x\|_1.$$
\end{itemize}

An important difference between classical and non-commutative
martingales is that the decomposition stated in Theorem A requires
four martingales versus the three martingales of Gundy's classical
decomposition. This difference is highlighted in Section
\ref{Section2} and is essentially due to the row and column nature
of Hardy spaces for non-commutative martingales from \cite{PX1}.
We also remark that Gundy's paper \cite{G} is based mainly on
stopping time arguments which at the time of this writing do not
appear to have a trackable non-commutative extension. Our approach
is based on a non-commutative analogue of Doob's maximal
inequality formulated by Cuculescu for positive martingales in
\cite{C}. Let us also mention that a weaker version of Gundy's
decomposition was obtained by Burkholder in \cite{B}. This
alternative decomposition $f=a+b+c$ does not satisfies the
$L_{\infty}$-estimate for the martingale $a$. However, it only
uses one stopping time (Gundy's approach uses two) in the
construction and is therefore easier to handle for many
applications. We shall also obtain a non-commutative analogue of
Burkholder's decomposition which will be used in some of the
applications we present in this paper.

\vskip5pt

As in the classical context, our decomposition is a powerful tool
to prove weak type inequalities for non-commutative martingales.
We illustrate this by reproving the main results in \cite{R} and
\cite{R3} respectively. More concretely, the weak type $(1,1)$
boundedness of non-commutative martingale transforms and the
non-commutative analogue of Burkholder's weak type inequality for
square functions. The latter result was recently proved in
\cite{R3} and can be regarded as the weak type extension of
non-commutative Burkholder-Gundy inequality from \cite{PX1}. The
contribution of our approach lies in the simplicity of the proofs,
derived from the new insight provided by Gundy's decomposition.

\vskip5pt

The last application is a non-commutative extension of a classical
result of Aldous and Fremlin \cite{ALFR} on basic sequences on
$L_1$-spaces. Recall that a basic sequence $(x_n)_{n \ge 1}$ in a
Banach space $\mathrm{X}$ is said to be $2$-co-lacunary  if there
is a constant $\delta>0$ so that for any finite sequence
$(a_n)_{n\geq 1}$ of scalars, $$\delta \Big(\sum_{n \ge 1} |a_n|^2
\Big)^{1/2} \le \Big\| \sum_{n \ge 1} a_n x_n
\Big\|_{\mathrm{X}}.$$ As application of Theorem A, we shall prove
that any relatively weakly compact martingale difference sequence
in $L_1(\M,\T)$ whose sequence of norms is bounded away from zero
is a $2$-co-lacunary sequence. Using this, we shall also prove
that for any semifinite and hyperfinite von Neumann algebra
$\mathcal{M}$, every bounded sequence in $L_1(\mathcal{M},\tau)$
has either a convergent or a 2-co-lacunary subsequence.

\vskip5pt

The paper is organized as follows. In Section \ref{Section1}, we
set some basic preliminary background concerning non-commutative
symmetric spaces and martingale theory that will be needed
throughout. Section~2 is devoted mainly to the statement and proof
of the main result along with some reformulations. In Section
\ref{Section3}, we present the three applications mentioned above.
Our notation and terminology are standard as may be found in the
books \cite{KR1} and \cite{T}. The letter $\mathrm{c}$ will denote
an absolute constant which might change from one instance to
another.

\vskip5pt

\noindent \textsc{Acknowledgment}. The first-named author would
like to express his gratitude to T. Mart\'{\i}nez and J.L. Torrea for
having brought Gundy's paper to his attention.

\section{Preliminary definitions and results}
\label{Section1}

This section is devoted to some preliminary definitions and
results that might be well-known to experts in the field and that
will be needed in the rest of the paper. Throughout, $\M$ is a
semifinite von Neumann algebra with a normal faithful semifinite
trace  $\T$. The identity element of $\M$ is denoted by ${\bf 1}$.
For $0 < p \leq \infty$, let $L_p(\M,\T)$ be the associated
non-commutative $L_p$-space, see for instance \cite{N,PX2}. Note
that if $p=\infty$, $L_\infty(\M,\tau)$ is just $\M$ with the
usual operator norm; also recall that for $0< p<\infty$, the
(quasi)-norm on $L_p(\M,\T)$ is defined by $$\Vert x \Vert_p =
\big( \T(|x|^p) \big)^{1/p}, \quad \mbox{where} \quad
|x|=(x^*x)^{1/2}.$$

\subsection{Non-commutative symmetric spaces}
Assume that $\M$ is acting on a Hilbert space $H$. A closed
densely defined operator $x$ on $H$ is \emph{affiliated} with
$\mathcal{M}$ if $x$ commutes with every unitary $u$ in the
commutant $\mathcal{M}'$ of $\cal{M}$. If $a$ is a densely defined
self-adjoint operator on $H$ and $a = \int_{\R} s d e^a_s$ is its
spectral decomposition, then for any Borel subset $B \subseteq
\R$, we denote by $\chi_B(a)$ the corresponding spectral
projection $\int_{\R} \chi_B(s) d e^a_s$. An operator $x$
affiliated with $\mathcal{M}$ is called \emph{$\tau$-measurable}
if there exists $s > 0$ such that $$\tau \big( \chi_{(s,\infty)}
(|x|) \big) < \infty.$$ The \emph{generalized singular-value}
$\mu(x): \R_+ \to \R_+$ of a $\tau$-measurable $x$ is defined by
$$\mu_t (x) = \inf \Big\{ s
> 0 \, \big| \ \tau \big( \chi_{(s,\infty)}(|x|) \big) \le t \Big\}.$$
The reader is referred to \cite{FK} for a  detailed exposition of
the function $\mu(x)$.  For  a rearrangement invariant
quasi-Banach function space $E$ on the interval $(0,
\tau(\mathbf{1}))$, we define the \emph{non-commutative symmetric
space} $E(\mathcal{M}, \tau)$ by setting: $$E(\mathcal{M},\tau) :=
\Big\{x \in L_0(\mathcal{M}, \tau) \, \big| \ \mu(x) \in E \Big\},
\quad \text{and} \quad \left\| x \right\|_{E(\M,\T)} :=
\left\|\mu(x)\right\|_E$$ where $L_0(\mathcal{M},\tau)$ stands for
the $*$-algebra of $\tau$-measurable operators. It is known that
$(E(\mathcal{M}, \tau), \|\cdot\|_{E(\M,\T)})$ is a Banach
(respectively, quasi-Banach) space whenever $E$ is a Banach
(respectively, quasi-Banach) function space. We refer the reader
to \cite{DDdP,X2} for more in depth discussion of this
construction. For the case where $E$ is a  Banach space, the
inclusions \[ L_1(\M,\T) \cap \M \subseteq E(\M,\T) \subseteq
L_1(\M,\T) + \M \] hold with the inclusion maps being of norm one
(here, the norms in $L_1(\M,\T) \cap \M$ and $L_1(\M,\T) + \M$ are
the usual norms of the intersection and sum of Banach spaces). The
\emph{K\"othe dual} $E(\M,\T)^\times$ of $E(\M,\T)$ is defined to
be the set of all $x \in L_0({\M},\T)$ such that $xy \in
L_1(\M,\T)$ for all $y \in E(\M,\T)$. With the norm defined by
setting: $$\left\|x\right\|_{E(\M,\T)^\times} := \sup \Big\{
\T(|xy|) \, \big| \ y \in E(\M,\T), \, \|y\|_{E(\M,\T)} \leq 1
\Big\},$$ the K\"othe dual $E(\M,\T)^\times$ is a Banach space.
Basic properties of K\"othe duality for the commutative case may
be found in \cite{LT}. For the non-commutative setting, the reader
is referred to  \cite{DDP3}. We note from \cite{DDP3} that if $E$
is a rearrangement invariant function space $E$ then
$(E(\M,\T)^\times, \|\cdot\|_{E(\M,\T)^\times})$ may be identified
with the space $(E^\times(\M,\T), \|\cdot\|_{E^\times(\M,\T)})$.
In particular,
\begin{eqnarray*}
(L_1(\M,\T) +\M)^\times & = & L_1(\M,\T) \cap \M, \\ (L_1(\M,\T)
\cap \M)^\times & = & L_1(\M,\T) +\M.
\end{eqnarray*}

Relative weak compactness in non-commutative spaces plays a role
in this paper. Below, we  explicitly state a  characterization
that we need in the subsequent sections. First, we set $S_0(\M,\T)
:={\M}_0 \cap (L_1(\M,\T) +\M)$ with
\[ {\M}_0 := \Big\{ x \in L_0({\M},\T) \, \big| \ \mu_t(x)\to 0 \
\text{as} \ t \to \infty \Big\}. \]

\begin{theorem} \label{weakcom2} {\rm\cite[Theorem~5.4]{DSS}.}
Assume that the symmetric space $E(\M,\T)$ is contained in
$S_0(\M,\T)$ and that $K$ is a bounded subset of
$E(\M,\T)^\times$. Then, the following statements are equivalent:
\begin{itemize}
\item[(i)] $\mu(K)$ is relatively $\sigma(E^\times, E)$-compact;
\item[(ii)] $K$ is relatively $\sigma(E(\M,\T)^\times,
E(\M,\T))$-compact.
\end{itemize}
\end{theorem}

Our interest in this paper is mainly restricted to non-commutative
$L_p$-spaces and \emph{non-commutative weak $L_1$-spaces}.
Following the construction of symmetric spaces of measurable
operators, the non-commutative weak $L_1$-space
$L_{1,\infty}(\mathcal{M}, \T )$,  is defined as the set of all
$x$ in $L_0(\mathcal{M},\tau)$ for which the quasi-norm
\[ \left\|x\right\|_{1,\infty} = \sup_{t > 0} t \mu_t(x) =
\sup_{\lambda > 0} \lambda \tau \big( \chi_{(\lambda, \infty)}
(|x|) \big) < \infty. \] As in the commutative case, it can be
easily verified that if $x_1, x_2 \in L_{1,\infty}(\M,\T)$ then
$\|x_1 +x_2\|_{1,\infty} \leq 2\|x_1\|_{1,\infty}
+2\|x_2\|_{1,\infty}$. In fact, the following more general
quasi-triangle inequality holds and will be used repeatedly in the
sequel. A short proof can be found in \cite[Lemma~1.2]{R3}.

\begin{lemma} \label{Quasi-Triangle}
Given two operators $x_1, x_2$ in $L_{1,\infty}(\mathcal{M},\T)$
and $\lambda > 0$, we have $$\lambda \, \tau \Big(
\chi_{(\lambda,\infty)} \big( |x_1+x_2| \big) \Big) \le 2 \lambda
\, \tau \Big( \chi_{(\lambda/2, \infty)} \big( |x_1| \big) \Big) +
2 \lambda \, \tau \Big( \chi_{(\lambda/2, \infty)} \big( |x_2|
\big) \Big).$$
\end{lemma}

Let $\mathsf{P} = \{ p_i \}_{i=1}^m$ be a finite sequence of
mutually orthogonal projections in $\mathcal{M}$. We consider the
{\it triangular truncation} with respect to $\mathsf{P}$ as the
mapping on $L_0(\mathcal{M},\tau)$ defined by:
\[ \mathcal{T}^{(\mathsf{P})} x = \sum_{i=1}^m \sum_{i\leq j} p_i x
p_j. \] The following  lemma will be used in the sequel.

\begin{lemma} \label{Truncation} {\rm\cite[Proposition~1.4]{R3}} There
exists an absolute constant $\mathrm{c}>0$ so that if
$(\mathsf{P}_k)_{k\geq 1}$ is a family of finite sequences of
mutually orthogonal projections and $(x_k)_{k\geq 1}$  is a
sequence in $L_1(\mathcal{M},\T)$, then
\[
 \Big\| \Big( \sum_{k\geq 1} \big| \mathcal{T}^{(\mathsf{P}_k)} x_k
\big|^2 \Big)^{1/2} \Big\|_{1,\infty} \le \mathrm{c} \sum_{k\geq
1} \|x_k\|_1.
\]
\end{lemma}

\subsection{Non-commutative martingales}

Consider a von Neumann subalgebra $\mathcal{N}$ of $\mathcal{M}$
(i.e. a weak$^*$ closed $*$-subalgebra of $\mathcal{M}$). A
\emph{conditional expectation} $\mathsf{E}: \mathcal{M} \to
\mathcal{N}$ from $\mathcal{M}$ onto $\mathcal{N}$ is a positive
contractive projection. The conditional expectation $\mathsf{E}$
is called \emph{normal} if the adjoint map $\mathsf{E}^*$
satisfies $\mathsf{E}^*(\mathcal{M}_*) \subset \mathcal{N}_*$. In
this case, there is map $\mathsf{E}_*: \mathcal{M}_* \rightarrow
\mathcal{N}_*$ whose adjoint is $\mathsf{E}$. Note that such
normal conditional expectation exists if and only if the
restriction of $\tau$ to the von Neumann subalgebra $\mathcal{N}$
remains semifinite (see for instance \cite[Theorem~3.4]{T}). Any
such conditional expectation is trace preserving (that is, $\tau
\circ \mathsf{E} = \tau$) and satisfies the bimodule property:
\[ \mathsf{E}(axb) = a \mathsf{E}(x) b \quad \mbox{for all} \quad a,b
\in \mathcal{N} \ \mbox{and} \ x \in \mathcal{M}.\]

Let $(\mathcal{M}_n)_{n\ge 1}$  be an increasing sequence of von
Neumann subalgebras of $\mathcal{M}$ such that the union of the
$\mathcal{M}_n$'s is weak$^*$ dense in $\mathcal{M}$. Assume that
for every $n\geq 1$, there is a normal conditional expectation
$\mathsf{E}_n: \mathcal{M} \to \mathcal{M}_n$. Note that for every
$1 \leq p<\infty$ and $n\ge 1$, $\mathsf{E}_n$ extends to a
positive contraction $\mathsf{E}_n: L_p(\mathcal{M},\T) \to
L_p(\mathcal{M}_n, \T|_{\M_n})$. A \emph{non-commutative
martingale} with respect to the filtration $(\mathcal{M}_n)_{n \ge
1}$ is a sequence $x = (x_n)_{n \ge 1}$ in $L_1(\mathcal{M},\T)$
such that
\[
\mathsf{E}_m(x_n) = x_m \quad \mbox{for all} \quad 1
\le m \le n < \infty.
\]
If additionally $x \subset L_p(\mathcal{M},\T)$  for some $1\leq
p\leq \infty$, then $x$ is called a \emph{$L_p$-martingale}. In
this case, we set
\[\left\|x\right\|_p:=\sup_{n \ge 1} \left\|x_n\right\|_p.
\]
If $\|x\|_p<\infty$, $x$ is called a \emph{$L_p$-bounded
martingale}. Given a martingale $x = (x_n)_{n \ge 1}$, we assume
the convention that $x_0 = 0$. Then, the martingale difference
sequence $dx =(dx_k)_{k \ge 1}$ associated to $x$ is defined by
$$dx_k = x_k - x_{k-1}.$$


We now describe square functions of non-commutative martingales.
Following Pisier and Xu \cite{PX1}, we will consider the following
row and column versions of square functions. Given a martingale
difference sequence $dx = (dx_k)_{k \ge 1}$ and  $n \geq 1$, we
define the \emph{row} and \emph{column square functions} of $x$ as
\[
\cal{S}_{C,n} (x) := \Big(\sum^n_{k=1}|dx_{k}|^{2} \Big)^{{1}/{2}}
\quad \mathrm{and} \quad \cal{S}_{R,n}(x) :=
\Big(\sum^n_{k=1}|dx_{k}^*|^{2} \Big)^{{1}/{2}}.
\]
Let us consider a rearrangement invariant (quasi)-Banach function
space $E$ on the interval $[0, \T({\bf 1}))$. Then, we define the
spaces $E(\mathcal{M}, \tau; l^2_C)$ and $E(\mathcal{M}, \tau;
l^2_R)$ as the completions of the vector space of finite sequences
$a = (a_k)_{k \ge 1}$ in $E(\mathcal{M}, \tau)$ with respect to
the following norms
\begin{eqnarray*}
\|a\|_{E(\M,\tau;l^{2}_{C})} & = & \Big\| \Big( \sum_{k \geq 1}
|a_{k}|^{2} \Big)^{{1}/{2}} \Big\|_{E(\M,\T)}, \\ \|a\|_{E(\M,
\tau; l^{2}_{R})} & = & \Big\| \Big( \sum_{k \geq 1}
|a_{k}^{*}|^{2} \Big)^{{1}/{2}} \Big\|_{E(\M,\T)}.
\end{eqnarray*}

The martingale difference sequence $dx$ belongs to
$E(\M,\tau;l^{2}_{C})$ (respectively,  $E(\M,\tau; l^{2}_{R})$) if
and only if the sequence $(\mathcal{S}_{C,n}(x))_{n \ge 1}$
(respectively, $(\mathcal{S}_{R,n}(x))_{n \ge 1}$) is bounded in
$E(\M,\T)$. In this case,  the limits
\[
\cal{S}_{C}(x) := \Big( \sum^{\infty}_{k=1} |dx_{k}|^{2}
\Big)^{{1}/{2}} \quad \mbox{and} \quad \mathcal{S}_{R}(x) := \Big(
\sum^{\infty}_{k=1}|dx_{k}^{*}|^{2} \Big)^{{1}/{2}} \] are
elements of $E(\M,\T)$. These two versions of square functions are
very crucial in the subsequent sections.

\vskip5pt


The next result is by now well known for positive martingales
\cite{C} and can be viewed as a non-commutative analogue of the
classical weak type $(1,1)$ boundedness of Doob's maximal
function. The extension to self-adjoint martingales stated below
plays a crucial role in the next section. Its proof is a minor
adjustment of the original argument of Cuculescu  but we will
include the details for completeness.

\begin{proposition}\label{Cuculescu}
If $x=(x_{n})_{n \ge 1}$ is a self-adjoint $L_1$-bounded
martingale and $\lambda$ is a positive real number, there exists a
sequence of decreasing projections in the von Neumann algebra $\M$
$$q_0^{(\lambda)} \ge q_1^{(\lambda)} \ge q_2^{(\lambda)} \ge
\ldots$$ satisfying the following properties:
\begin{itemize}
\item[(i)] for every $n \geq 1$, $q_{n}^{(\lambda)} \in \M_{n}$;
\item[(ii)] for every $n \geq 1$, $q_{n}^{(\lambda)}$ commutes
with $q_{n-1}^{(\lambda)} x_{n} q_{n-1}^{(\lambda)}$;
\item[(iii)] for every $n \geq 1$, $|q_{n}^{(\lambda)} x_{n}
q_{n}^{(\lambda)}| \leq \lambda q_{n}^{(\lambda)}$;
\item[(iv)]
if we  set $q^{(\lambda)} = \bigwedge_{n=1}^{\infty}
q_{n}^{(\lambda)}$, then $$\T \Big( {\bf 1}-q^{(\lambda)} \Big)
\leq \frac{1}{\lambda} \|x\|_1.$$
\end{itemize}
\end{proposition}

\dem Let $q_{0}^{(\lambda)} = {\bf 1}$ and inductively on $n \geq
1$, define
\[
q_{n}^{(\lambda)} := q_{n-1}^{(\lambda)} \chi_{[-\lambda,\lambda]}
\Big( q_{n-1}^{(\lambda)} x_{n} q_{n-1}^{(\lambda)} \Big) =
\chi_{[-\lambda,\lambda]} \Big( q_{n-1}^{(\lambda)} x_{n}
q_{n-1}^{(\lambda)} \Big) q_{n-1}^{(\lambda)}.\]
The identity above follows since $q_{n-1}^{(\lambda)}$ commutes
with $q_{n-1}^{(\lambda)} x_{n} q_{n-1}^{(\lambda)}$ because
$q_{n-1}^{(\lambda)}$ is a projection. This clearly gives a
decreasing sequence of projections. By induction, condition (i)
holds. Moreover, condition (ii) follows directly from the
definition above. For (iii), note that for every $n\geq 1$,
\begin{eqnarray*}
q_{n}^{(\lambda)} x_{n} q_{n}^{(\lambda)} & = & q_{n}^{(\lambda)}
(q_{n-1}^{(\lambda)} x_{n} q_{n-1}^{(\lambda)}) q_{n}^{(\lambda)}
\\ & = & q_{n-1}^{(\lambda)} \chi_{[-\lambda,\lambda]}
(q_{n-1}^{(\lambda)} x_{n} q_{n-1}^{(\lambda)})
q_{n-1}^{(\lambda)} x_{n} q_{n-1}^{(\lambda)}
\chi_{[-\lambda,\lambda]} (q_{n-1}^{(\lambda)} x_{n}
q_{n-1}^{(\lambda)}) q_{n-1}^{(\lambda)}.
\end{eqnarray*}
Therefore $-\lambda q_{n}^{(\lambda)}\leq q_{n}^{(\lambda)} x_{n}
q_{n}^{(\lambda)}\leq \lambda q_{n}^{(\lambda)}$ and (iii)
follows. To prove (iv), we use the non-commutative analogue of
Krickeberg's decomposition \cite{C}, so that we may write $x_n=
w_n - z_n$ where $w = (w_n)_{n \ge 1}$ and $z = (z_n)_{n \ge 1}$
are positive martingales with
$$\|x\|_1=\T(w_1 + z_1).$$ For every $n\geq 1$,
\begin{eqnarray*}
\|x\|_1 & = & \T \big( (w_n + z_n) q_{n}^{(\lambda)} \big) +
\sum^n_{k=1} \T \Big( (w_n + z_n) (q_{k-1}^{(\lambda)} -
q_{k}^{(\lambda)}) \Big)
\\ & = & \T \big( q_{n}^{(\lambda)} (w_n + z_n)
q_n^{(\lambda)} \big) + \sum^n_{k=1} \T \Big( \mathsf{E}_k (w_n +
z_n) (q_{k-1}^{(\lambda)} - q_{k}^{(\lambda)}) \Big).
\end{eqnarray*}
Since $\T \big( q_{n}^{(\lambda)} (w_n + z_n) q_{n}^{(\lambda)})
\geq 0$, we have
\begin{eqnarray*}
\|x\|_1 & \geq & \sum^{n}_{k=1} \T \Big( (q_{k-1}^{(\lambda)} -
 q_{k}^{(\lambda)}) (w_k + z_k) (q_{k-1}^{(\lambda)} -
 q_{k}^{(\lambda)}) \Big) \\ & \geq & \T \Big( \sum^{n}_{k=1}
 \Big| (q_{k-1}^{(\lambda)} - q_{k}^{(\lambda)}) (w_k - z_k)
 (q_{k-1}^{(\lambda)} - q_{k}^{(\lambda)}) \Big| \Big)
\\ & = & \T \Big( \sum^{n}_{k=1} \Big| (q_{k-1}^{(\lambda)} -
q_{k}^{(\lambda)}) (q_{k-1}^{(\lambda)} x_k q_{k-1}^{(\lambda)})
(q_{k-1}^{(\lambda)} - q_{k}^{(\lambda)}) \Big| \Big).
\end{eqnarray*}
From the definition of $q_{k}^{(\lambda)}$, it is clear that
\begin{eqnarray*}
q_{k-1}^{(\lambda)} - q_{k}^{(\lambda)} & = & q_{k-1}^{(\lambda)}
\Big( \chi_{(-\infty, -\lambda)} (q_{k-1}^{(\lambda)} x_{k}
q_{k-1}^{(\lambda)}) + \chi_{(\lambda, \infty)}
(q_{k-1}^{(\lambda)} x_{k} q_{k-1}^{(\lambda)}) \Big) \\ & = &
\Big( \chi_{(-\infty, -\lambda)} (q_{k-1}^{(\lambda)} x_{k}
q_{k-1}^{(\lambda)}) + \chi_{(\lambda, \infty)}
(q_{k-1}^{(\lambda)} x_{k} q_{k-1}^{(\lambda)}) \Big)
q_{k-1}^{(\lambda)}.
\end{eqnarray*}
Therefore, if $q_{k-1}^{(\lambda)}
x_{k}q_{k-1}^{(\lambda)}=\int_{\R} t de_t^{(k)}$ is the spectral
decomposition of $q_{k-1}^{(\lambda)} x_{k}q_{k-1}^{(\lambda)}$,
we find that
\begin{eqnarray*}
\lambda (q_{k-1}^{(\lambda)} - q_{k}^{(\lambda)}) & \le &
\int_{-\infty}^{- \lambda} |t| \, de_t^{(k)} +
\int_{\lambda}^{\infty} |t| \, de_t^{(k)} \\ & = & \Big|
(q_{k-1}^{(\lambda)} - q_{k}^{(\lambda)}) (q_{k-1}^{(\lambda)} x_k
q_{k-1}^{(\lambda)}) (q_{k-1}^{(\lambda)} - q_{k}^{(\lambda)})
\Big|.
\end{eqnarray*}
We can now conclude that $$\T \Big( {\bf 1}-q_{n}^{(\lambda)}
\Big) \le \frac{1}{\lambda} \|x\|_1.$$ Taking the limit as $n \to
\infty$, we obtain (iv). This completes the proof. \fin

In the following, we will refer to the sequence of projections of
Proposition \ref{Cuculescu} as the sequence of Cuculescu's
projections associated to the (self-adjoint) martingale $x$ and
the (positive) parameter $\lambda$. In the next result, we collect
some basic properties of this sequence that are very useful for
the presentation in the next section.

\begin{proposition} \label{Randri}
Let $x = (x_n)_{n \ge 1}$ be a  self-adjoint $L_1$-bounded
martingale and $\lambda$ a positive real number. Then, the
sequence of Cuculescu's projections associated to $x$ and
$\lambda$ satisfies the following estimates for every $n\ge 1$:
\begin{eqnarray}
\label{Est1} \sum_{k=1}^{n} \big\| q_{k-1}^{(\lambda)} x_k
q_{k-1}^{(\lambda)} - q_k^{(\lambda)} x_k q_k^{(\lambda)} \big\|_1
& \le & \|x\|_1, \\ \label{Est2} \sum_{k=1}^{n} \big\|
q_{k-1}^{(\lambda)} x_{k-1} q_{k-1}^{(\lambda)} - q_k^{(\lambda)}
x_{k-1} q_k^{(\lambda)}\big\|_1 & \le & 2 \|x\|_1, \\ \label{Est3}
\sum_{k=1}^{n} \big\| q_{k-1}^{(\lambda)} dx_k q_{k-1}^{(\lambda)}
- q_k^{(\lambda)} dx_k q_k^{(\lambda)} \big\|_1 & \le & 3 \|x\|_1.
\end{eqnarray}
Moreover, the following identity holds
\begin{equation} \label{Ident1}
\sum_{k=1}^{n} q_{k-1}^{(\lambda)} dx_k  q_{k-1}^{(\lambda)}=
q_n^{(\lambda)}  x_n  q_n^{(\lambda)}+ \sum_{k=1}^{n}
(q_{k-1}^{(\lambda)}-q_{k}^{(\lambda)})x_k
(q_{k-1}^{(\lambda)}-q_{k}^{(\lambda)}).
\end{equation}
In particular, we obtain
\begin{equation} \label{Est4}
\Big\| \sum_{k=1}^{n} q_{k-1}^{(\lambda)} dx_k q_{k-1}^{(\lambda)}
\Big\|_1 \leq 2 \|x\|_1.
\end{equation}
\end{proposition}

\dem We will write $(q_n)_{n \ge 0}$ for $(q_n^{(\lambda)})_{n \ge
0}$ and $q$ for $q^{(\lambda)}$ (see Proposition \ref{Cuculescu}).
Let $v_k = q_{k-1} x_k q_{k-1} - q_k x_k q_k$. Since $q_k$
commutes with $q_{k-1} x_k q_{k-1}$, we have $v_k = (q_{k-1} -
q_k) x_k (q_{k-1} - q_k)$ and therefore for any given $n \ge 1$,
\begin{eqnarray*}
\sum_{k=1}^n \|v_k\|_1 & = & \sum_{k=1}^n \big\| (q_{k-1} - q_k)
x_k (q_{k-1} - q_k) \big\|_1 \\ & = & \sum_{k=1}^n \big\|
\mathsf{E}_k \big( (q_{k-1} - q_k) x_n (q_{k-1} - q_k) \big)
\big\|_1
 \\ & \leq & \sum_{k=1}^n \big\| (q_{k-1} - q_k) x_n (q_{k-1}
- q_k) \big\|_1 \le \|x_n\|_1.
\end{eqnarray*}
Hence, inequality (\ref{Est1}) is satisfied.  For inequality
(\ref{Est2}), set $$\sigma_k = q_k x_{k-1} q_k - q_{k-1} x_{k-1}
q_{k-1}.$$ Then we have,
\begin{eqnarray*}
\sigma_k & = & q_k x_{k-1} (q_k - q_{k-1}) + (q_k -q_{k-1})
x_{k-1} q_{k-1} \\ & = & q_k q_{k-1} x_{k-1} q_{k-1} (q_k -
q_{k-1}) + (q_k - q_{k-1}) q_{k-1} x_{k-1} q_{k-1}.
\end{eqnarray*}
By H\"{o}lder's inequality and Proposition~\ref{Cuculescu}, we
deduce that
\[\sum_{k=1}^\infty\|\sigma_k\|_1 \le 2 \sum_{k=1}^\infty
\tau(q_{k-1} - q_k) \big\| q_{k-1} x_{k-1} q_{k-1}\big\|_{\infty}
\le 2 \lambda \tau (\mathbf{1} - q)\leq 2\|x\|_1,
 \]
which proves inequality (\ref{Est2}). The estimate in (\ref{Est3})
follows directly from (\ref{Est1}) and (\ref{Est2}) and the
triangle inequality. The identity (\ref{Ident1})  follows
immediately from summing by parts. Indeed, for  $n\geq 1$ we have
\begin{eqnarray*}
\sum_{k=1}^n q_{k-1} dx_k q_{k-1} & = & \sum_{k=1}^n \big( q_{k-1}
x_k q_{k-1} - q_{k-1} x_{k-1} q_{k-1} \big) \\ & = &
\sum_{k=1}^{n-1} \big( q_{k-1} x_k q_{k-1} - q_k x_k q_k \big) +
q_{n-1} x_n q_{n-1} \\ & = & \sum_{k=1}^{n} \big( q_{k-1} x_k
q_{k-1} - q_k x_k q_k \big) + q_{n} x_n q_{n}
\end{eqnarray*}
Finally, (\ref{Est4}) follows from (\ref{Est1}) and
(\ref{Ident1}). The proof is complete. \fin

\section{Non-commutative Gundy's decomposition}
\label{Section2}

In this section we present the non-commutative analogue of Gundy's
theorem, which is the main result of this paper. All adapted
sequences and martingales are understood to be with respect to a
fixed filtration $(\mathcal{M}_n)_{n \ge 1}$ of von Neumann
subalgebras of $\M$. For convenience, we assume that
$\mathsf{E}_0=\mathsf{E}_1$.

\begin{theorem}\label{Main} If $x=(x_n)_{n \ge 1}$ is a $L_1$-bounded
non-commutative martingale and $\lambda$ is a positive real
number, there exist four martingales $\alpha$, $\beta$, $\gamma$,
and $\upsilon$ satisfying the following properties for some
absolute constant $\mathrm{c}$:
\begin{itemize}
\item[(i)] $x=\alpha +\beta + \gamma + \upsilon$;
\item[(ii)] the martingale $\alpha$ satisfies $$\|\alpha\|_1 \leq
\mathrm{c} \|x\|_1, \quad \|\alpha\|_2^2 \leq \mathrm{c}
\lambda\|x\|_1, \quad \|\alpha\|_\infty \leq \mathrm{c} \lambda;$$
\item[(iii)] the martingale $\beta$ satisfies $$\sum_{k=1}^{\infty}
\|d\beta_k\|_1 \le \mathrm{c} \|x\|_1;$$
\item[(iv)] $\gamma$ and $\upsilon$ are $L_1$-martingales with
$$\max \Big\{ \lambda \tau \Big( \bigvee_{k \ge 1} \mathrm{supp}
|d\gamma_k| \Big), \, \lambda \tau \Big( \bigvee_{k \ge 1}
\mathrm{supp} \, |d\upsilon_k^*| \Big) \Big\} \le \mathrm{c}
\|x\|_1.$$
\end{itemize}
\end{theorem}

\dem Without loss of generality, we can and do  assume that the
martingale $x$ is positive. Denote by $(q_n)_{n \ge 0}$ the
sequence of Cuculescu's projections associated with the martingale
$x$ and a fixed $\lambda>0$. The construction is done in two
steps.

\vskip5pt

\noindent{\bf Step~1.} We consider the following martingale
difference sequence
$$dy_k := q_k dx_k q_k - \E_{k-1} (q_k dx_k q_k ) \quad \text{for}
\quad k \geq 1.$$
It is clear that $(dy_k)_{k \ge 1}$ is a martingale difference
sequence and the corresponding martingale $y = (y_n)_{n \ge 1}$ is
a self-adjoint $L_1$-martingale. The following intermediate lemma
is essential for our construction.
\begin{lemma}\label{marting-y}
The martingale $y$ is $L_1$-bounded with $\|y\|_1\leq  9 \|x\|_1$.
\end{lemma}
\dem This follows essentially from Proposition~\ref{Randri}.
Indeed, for every $n\geq 1$,
\begin{eqnarray*}
\left\|y_n\right\|_1 & = & \Big\| \sum_{k=1}^n q_k dx_k q_k -
\E_{k-1} (q_k dx_k q_k ) \Big\|_1 \\ & \leq & \Big\|\sum_{k=1}^n
q_{k-1} dx_k q_{k-1} \Big\|_1 \\ & + & \sum_{k=1}^n \big\| q_{k-1}
dx_k q_{k-1}-q_k dx_k q_k \big\|_1  \\ & + & \sum_{k=2}^n \big\|
\E_{k-1} (q_k dx_k q_k) \big\|_1 + \big\| q_1 x_1 q_1 \big\|_1.
\end{eqnarray*}
Since for $2\leq k\leq n$, $$\mathsf{E}_{k-1} \big( q_k dx_k q_k
\big) = \mathsf{E}_{k-1} \big( q_k dx_k q_k - q_{k-1} dx_k q_{k-1}
\big),$$ the assertion follows from estimates (\ref{Est3}) and
(\ref{Est4}). Thus the lemma is proved. \fin

\noindent{\bf Step~2.} Let $(\pi_n)_{n \ge 0}$ stand for the
sequence of Cuculescu's projections relative to the (self-adjoint)
martingale $y$ and the parameter $\lambda$ fixed above. We define
the martingales $\alpha$, $\beta$, $\gamma$, and $\upsilon$ as
follows:

\begin{equation*}\label{decomposition}
\begin{cases}
d\alpha_k &:= \pi_{k-1} \big[ q_k dx_k q_k - \E_{k-1} (q_k dx_k
q_k) \big] \pi_{k-1}, \\ d\beta_k &:= \pi_{k-1} \big[ q_{k-1}dx_k
q_{k-1}- q_k dx_k q_k + \E_{k-1} (q_k dx_k q_k) \big]\pi_{k-1}, \\
d\gamma_k &:=dx_k -  dx_k q_{k-1}\pi_{k-1},\\ d\upsilon_k &:=dx_k
q_{k-1}\pi_{k-1}-\pi_{k-1}q_{k-1}dx_k q_{k-1}\pi_{k-1}.
\end{cases} \tag{$\mathbf{G}_{\lambda}$}
\end{equation*}
Clearly, $d\alpha$, $d\beta$, $d\gamma$ and $d\upsilon$ are
martingale difference sequences and $x=\alpha +\beta +\gamma
+\upsilon$.

\begin{lemma}\label{alpha-norm} The martingale $\alpha$ satisfies
$$\|\alpha\|_1 \leq 18 \|x\|_1, \quad \|\alpha\|_2^2 \leq 72
\lambda\|x\|_1, \quad \|\alpha\|_\infty \leq 4\lambda.$$
\end{lemma}

\dem Note that for $k \geq 1$,  $d\alpha_k =
\pi_{k-1}dy_k\pi_{k-1}$. Thus, the $L_1$-estimate follows directly
from (\ref{Est4}) and Lemma~\ref{marting-y}. For the
$L_\infty$-estimate, we recall from (\ref{Ident1}) that for every
$n\geq 1$,
\[
\alpha_n= \sum_{k=1}^{n} \pi_{k-1} dy_k  \pi_{k-1} = \pi_n  y_n
\pi_n + \sum_{k=1}^{n} (\pi_{k-1}-\pi_{k}) y_k
(\pi_{k-1}-\pi_{k}).
 \]
The key observation is that $\sup_{k\geq 1} \|dy_k\|_\infty \leq 2
\lambda$. We have the following estimates:
\begin{eqnarray*}
\left\|\alpha_n \right\|_\infty & \leq & \left\|\pi_n  y_n \pi_n
\right\|_\infty \\ & + & \Big\| \sum_{k=1}^{n} (\pi_{k-1}-\pi_{k})
dy_{k} (\pi_{k-1}-\pi_{k}) \Big\|_\infty
\\ & + & \Big\|\sum_{k=1}^{n}
(\pi_{k-1}-\pi_{k}) y_{k-1} (\pi_{k-1}-\pi_{k}) \Big\|_\infty.
\end{eqnarray*}
We deduce from the definition of $(\pi_n)_{n \ge 0}$ that,
$$\left\|\alpha_n \right\|_\infty \leq \lambda + \sup_{k \leq n}
\left\|dy_{k}\right\|_\infty + \sup_{k\leq n}
\left\|(\pi_{k-1}-\pi_{k}) y_{k-1}
(\pi_{k-1}-\pi_{k})\right\|_\infty \leq 4\lambda.$$ The
$L_2$-estimate follows from those of $L_1$ and $L_\infty$ using
H\"older's inequality. \fin

\begin{lemma}\label{beta-norm}
The martingale $\beta$ satisfies $$\sum_{k=1}^\infty
\|d\beta_k\|_1 \leq 7\|x\|_1.$$
\end{lemma}

\dem From the definition of $(d\beta_k)_{k \ge 1}$, we have
\[
\sum_{k=1}^{\infty} \left\|d\beta_k\right\|_1 \le
\left\|d\beta_1\right\|_1 + \sum_{k=2}^{\infty} \big\| q_{k-1}
dx_k q_{k-1} - q_k dx_k q_k \big\|_1 + \sum_{k=2}^{\infty} \big\|
\mathsf{E}_{k-1} \big( q_k dx_k q_k \big) \big\|_1.
\]
Since for every $k\geq 2$, $$\mathsf{E}_{k-1} \big( q_k dx_k q_k
\big) = \mathsf{E}_{k-1} \big( q_k dx_k q_k - q_{k-1} dx_k q_{k-1}
\big),$$  we conclude from (\ref{Est3}) that
\[
\sum_{k=1}^{\infty} \left\|d\beta_k\right\|_1 \le
\left\|x_1\right\|_1 + 2 \sum_{k=2}^{\infty} \big\| q_{k-1} dx_k
q_{k-1} - q_k dx_k q_k \big\|_1 \le 7 \|x\|_1.
\]
Thus we have the estimate as stated. \fin

\begin{lemma}\label{gamma-support}
For every $k\geq 1$,
\begin{itemize}
\item[(a)]$\mathrm{supp} |d\gamma_k| \leq {\bf 1}- \pi_{k-1}
\wedge q_{k-1}$;
\item[(b)] $\mathrm{supp} |d\upsilon_k^*| \leq {\bf 1}-
\pi_{k-1} \wedge q_{k-1}$.
\end{itemize}
Consequently, the martingales $\gamma$ and $\upsilon$ satisfy
$$\max \Big\{ \lambda \tau \Big( \bigvee_{k \ge 1} \mathrm{supp}
\, |d\gamma_k| \Big), \, \lambda \tau \Big( \bigvee_{k \ge 1}
\mathrm{supp} \, |d\upsilon_k^*| \Big) \Big\} \le 10 \, \|x\|_1.$$
\end{lemma}

\dem It is immediate from (\ref{decomposition}) that
$d\gamma_k=d\gamma_k({\bf 1}-\pi_{k-1}\wedge q_{k-1})$ and using
polar decomposition, we obtain $|d\gamma_k|=|d\gamma_k|({\bf
1}-\pi_{k-1}\wedge q_{k-1})$ which shows that $\mathrm{supp} \,
|d\gamma_k| \leq {\bf 1}- \pi_{k-1} \wedge q_{k-1}$. As a
consequence,
\[
\bigvee_{k\geq 1} \mathrm{supp} \, |d\gamma_k| \leq {\bf 1}- \pi
\wedge q.\] Therefore, we deduce $$\T \Big( \bigvee_{k\geq 1}
\mathrm{supp} \, |d\gamma_k| \Big) \leq \T({\bf 1}- \pi) + \T({\bf
1}-q) \leq \frac{1}{\lambda} \big( \|y\|_1 + \|x\|_1 \big) \leq
\frac{10}{\lambda} \, \|x\|_1.$$ The same argument applies to the
martingale difference sequence $(d\upsilon^*_k)_{k \ge 1}$. \fin

It is now clear that  by combining Lemma~\ref{alpha-norm},
Lemma~\ref{beta-norm}, and Lemma~\ref{gamma-support}, all the
estimates from items (ii), (iii), and (iv) of Theorem~\ref{Main}
are verified. Moreover, the fact that $\gamma$ and $\upsilon$ are
$L_1$-martingales follows directly from (\ref{decomposition}).
\fin


\begin{remark}
\emph{It is important to note that  in strong contrast  with the
commutative case, the consideration of  the fourth martingale
$\upsilon$  is necessary in the decomposition stated in
Theorem~\ref{Main}. Indeed, assume that the decomposition in
Theorem~\ref{Main} can be done with only three martingales. That
is, if for every $L_1$-bounded martingale $x$ and $\lambda>0$,
there is a decomposition $x=\alpha +\beta +\gamma$ satisfying
(ii), (iii) and (iv) of Theorem~\ref{Main}. Then a straightforward
adjustment of the argument from the classical case used in
\cite{G} would prove that there is an absolute constant
$\mathrm{c}$ such that $$\max \Big\{
\|\mathcal{S}_R(x)\|_{1,\infty}, \,
\|\mathcal{S}_C(x)\|_{1,\infty} \Big\} \leq \mathrm{c} \,
\|x\|_1.$$ In particular, a standard use of the real interpolation
method shows that there is a constant $\mathrm{c}_p$ depending
only on $p$ such that, if $x$ is a $L_p$-bounded martingale with
$1 < p < 2$, then we have $$\max \Big\{ \|\mathcal{S}_R(x)\|_{p},
\, \|\mathcal{S}_C(x)\|_{p} \Big\} \leq \mathrm{c}_p \, \|x\|_p.$$
This is in a direct conflict with the non-commutative analogue of
Burkholder-Gundy inequality proved in \cite{PX1}. Hence in
general, decomposition into three martingales is not possible in
Theorem~\ref{Main}. This observation confirms the relation between
our decomposition and the row/column nature of Hardy spaces for
non-commutative martingales. Moreover, a detailed inspection of
the arguments sketched above shows that in fact, the martingales
$\gamma$ and $\upsilon$ can be regarded as the \lq column\rq{} and
\lq row\rq{} part of its commutative counterpart in \cite{G}.}
\end{remark}

\begin{remark}
\emph{In the construction (\ref{decomposition}) above, we have
$\alpha_1=0$, $\beta_1= x_1$, $\gamma_1=0$, and $\upsilon_1=0$.
Any other choice of these first terms  could have been taken
without any difference on the properties stated in
Theorem~\ref{Main}. Our choice is motivated in part by our second
application below, where we need to have $\gamma_1 = \upsilon_1 =
0$.}
\end{remark}

In the next formulation, we observe that if one wants to use three
martingales in the decomposition of Theorem~\ref{Main}, then we
have to consider a weaker notion of support projections.

\begin{definition}
For a non necessarily self-adjoint operator $x \in \M$, we define
the \emph{two-sided null projection} of $x$ to be the greatest
projection $q$ satisfying $qxq = 0$. In this case, we set
$\mbox{supp}^* x := {\bf 1} - q$.
\end{definition}

Clearly,  $\mbox{supp} \, x = \mbox{supp}^* x$ if $\mathcal{M}$ is
abelian. In general, $\mbox{supp}^*$ is smaller than the usual
support in the sense that $\mbox{supp}^* x \le \mbox{supp} \, x$
for any self-adjoint $x \in \M$  and for a non-self adjoint $x \in
\M$, $\mbox{supp}^* x$ is a subprojection of both the right and
left supports of $x$.  Using  this notion of support projections,
we can state:

\begin{corollary}
If $x=(x_n)_{n \ge 1}$ is a $L_1$-bounded non-commutative
martingale and $\lambda$ is a positive real number, there exist
three martingales $a$, $b$, and $c$ satisfying the following
properties for some absolute constant $\mathrm{c}$:
\begin{itemize}
\item[(i)] $x=a+b+c$;
\item[(ii)] the martingale $a$ satisfies $$\|a\|_1 \leq \mathrm{c}
\|x\|_1, \quad \|a\|_2^2 \leq \mathrm{c} \lambda \|x\|_1, \quad
\|a\|_\infty \leq \mathrm{c} \lambda;$$
\item[(iii)] the martingale $b$ satisfies $$\sum_{k=1}^{\infty}
\|db_k\|_1 \le \mathrm{c} \|x\|_1;$$
\item[(iv)] the martingale $c$ satisfies
$$\lambda \tau \Big( \bigvee_{k \ge 1} \mathrm{supp}^* dc_k \Big)
\le \mathrm{c} \|x\|_1.$$
\end{itemize}
\end{corollary}

\dem According to (\ref{decomposition}), it is enough to set $a
:=\alpha$, $b:=\beta$ and $c:=\gamma +\upsilon$. Then, (ii) and
(iii) follow directly from Theorem~\ref{Main}. For (iv), we note
from (\ref{decomposition}) that for every $k\geq 1$, $dc_k=dx_k
-\pi_{k-1}q_{k-1}dx_kq_{k-1}\pi_{k-1}$. Thus we deduce that for
$k\geq 1$,  $(\pi_{k-1} \wedge q_{k-1})dc_k(\pi_{k-1} \wedge
q_{k-1})=0$ and therefore, $\mathrm{supp}^* dc_k \leq {\bf
1}-(\pi_{k-1} \wedge q_{k-1})$. In particular, we obtain $$\T
\Big( \bigvee_{k \ge 1} \mathrm{supp}^* dc_k \Big) \leq \T({\bf
1}-\pi) +\T({\bf 1}-q).$$ At this point, (iv) follows as in
Lemma~\ref{gamma-support}. The proof is complete.\fin

Let us note that Gundy's original proof in \cite{G} uses two
stopping times. This essentially explains why we need two steps in
our construction of (\ref{decomposition}). In \cite{B}, Burkholder
provided a weaker version of Gundy's decomposition where only the
$L_2$-estimate is required for the first martingale in his
decomposition. His approach uses only one stoping time. In the
next result, we provide a non-commutative analogue of Burkholder's
approach. This provide us with a simpler decomposition which is
more useful for some applications.


\begin{corollary} \label{Gundy}
Let $x=(x_n)_{n \ge 1}$ be a $L_1$-bounded positive martingale and
$\lambda$ be a positive real number. Let us consider the
decomposition of $x$ as a sum of four martingales $x= \alpha' +
\beta' + \gamma' + \upsilon'$ with martingale differences given by
\begin{equation*}\label{Decomposition2}
\begin{cases}
 d\alpha'_k &:= q_k^{(\lambda)} dx_k q_k^{(\lambda)} -
\E_{k-1} \big( q_k^{(\lambda)} dx_k q_k^{(\lambda)} \big), \\
d\beta'_k &:= q_{k-1}^{(\lambda)} dx_k q_{k-1}^{(\lambda)} -
q_k^{(\lambda)} dx_k q_k^{(\lambda)} + \E_{k-1} \big(
q_k^{(\lambda)} dx_k q_k^{(\lambda)} \big), \\ d\gamma'_k &:= dx_k
- dx_k q_{k-1}^{(\lambda)},\\ d\upsilon'_k &:= dx_k
q_{k-1}^{(\lambda)} - q_{k-1}^{(\lambda)} dx_k
q_{k-1}^{(\lambda)}.
\end{cases} \tag{$\mathbf{G}'_{\lambda}$}
\end{equation*}
Then, the following properties hold:
\begin{itemize}
\item[(i)] the martingale $\alpha'$ satisfies $$\|\alpha'\|_1 \le
\mathrm{c} \|x\|_1, \quad \|\alpha'\|_2^2 \leq \mathrm{c}
\lambda\|x\|_1;$$
\item[(ii)] the martingales $\beta', \gamma'$ and $\upsilon'$
behave as $\beta, \gamma$ and $\upsilon$ in Theorem $\ref{Main}$.
\end{itemize}
\end{corollary}

\dem The estimates for $\beta', \gamma'$ and $\upsilon'$ can be
verified verbatim as in the proof of Theorem~\ref{Main}. For the
$L_1$-estimate of $\alpha'$, we use Lemma \ref{marting-y}. To
estimate the $L_2$-norm of $\alpha'$, we note by orthogonality
that for $n\geq 1$,
\begin{equation*}
\left\|\alpha'_n\right\|_2^2
=\sum_{k=1}^n\left\|d\alpha'_k\right\|_2^2 \leq 4\sum_{k=1}^n
\left\|q_kdx_kq_k\right\|_2^2.
\end{equation*}
On the other hand,  since for every $k\geq 1$,
$q_kdx_kq_k=q_k(q_kx_kq_k -q_{k-1}x_{k-1}q_{k-1})q_k$, we get that
\begin{equation*}
\left\|\alpha'_n\right\|_2^2 \leq 4 \sum_{k=1}^n \left\|q_k x_kq_k
-q_{k-1}x_{k-1}q_{k-1}\right\|_2^2.
\end{equation*}
Finally, according to \cite[Lemma~3.4]{R}, this gives
$\left\|\alpha'_n \right\|_2^2 \leq 24\lambda\left\|x\right\|_1.$
\fin

\begin{remark} \label{Non-positive}
\emph{Corollary \ref{Gundy} trivially extends to non-positive
martingales.}
\end{remark}

\section{Applications}
\label{Section3}

\subsection{Non-commutative martingale transforms}

As a first application of our decomposition, we provide a very
simple proof of the weak type $(1,1)$ boundedness of
non-commutative martingale transforms obtained in \cite{R}. This
has implications in non-commutative martingale theory as well as
for estimating $\mathrm{UMD}$ constants of certain non-commutative
function spaces. The reader is referred to \cite{R} and to Xu's
survey \cite{X} for a detailed exposition of these implications.

\begin{theorem}\label{Transform}
There exists an absolute constant $\mathrm{c}$ such that for every
martingale $x=(x_k)_{k \ge 1}$ bounded in $L_1(\mathcal{M}, \tau)$
and every sequence $(\xi_{k})_{k \ge 0}$ in $\M$ satisfying the
following properties:
\begin{itemize}
\item[(i)] $\xi_0={\bf 1}$;
\item[(ii)] $\sup_{k\geq 1}\|\xi_{k}\|_\infty\leq 1$;
\item[(iii)] $\xi_{k-1} \in \mathcal{M}_{k-1} \cap \M'_{k}$ for $k\ge
1$;
\end{itemize}
the following estimate holds for all $n \ge 1$,
\begin{equation*}
\Big\| \sum^{n}_{k=1} \xi_{k-1} dx_k \Big\|_{1,\infty} \leq
{\mathrm{c}} \|x\|_{1}.
\end{equation*}
\end{theorem}

\dem We have to show  that
\begin{equation}\label{transform-Ineq}
\lambda \, \tau \Big( \chi_{(\lambda,\infty)} \Big( \Big|
\sum_{k=1}^n \xi_{k-1} dx_k \Big| \Big)\Big) \le \mathrm{c}
\|x\|_1,
\end{equation}
for every $0 < \lambda < \infty$. For this, we fix $\lambda > 0$
and consider the  decomposition $x = \alpha + \beta + \gamma
+\upsilon$ of $x$ associated to $\lambda$ from Theorem~\ref{Main}.
Using the elementary inequality $|a+b|^2 \le 2 |a|^2 + 2 |b|^2$
for operators, we have
\begin{eqnarray*}
\Big|\sum_{k=1}^n \xi_{k-1} dx_k \Big|^2 & \le & 4
\Big|\sum_{k=1}^n \xi_{k-1} d\alpha_k \Big|^2 + 4
\Big|\sum_{k=1}^n \xi_{k-1} d\beta_k \Big|^2 \\ & + & 4 \Big|
\sum_{k=1}^n \xi_{k-1} d\gamma_k \Big|^2 + 4 \Big| \sum_{k=1}^n
\xi_{k-1} d\upsilon_k \Big|^2.
\end{eqnarray*}
Taking the trace, we obtain from Lemma \ref{Quasi-Triangle}
\begin{eqnarray*}
\lefteqn{\lambda \tau \Big( \chi_{(\lambda,\infty)} \Big( \Big|
\sum_{k=1}^n \xi_{k-1} dx_k \Big| \Big) \Big)=\lambda  \tau \Big(
\chi_{(\lambda^2,\infty)} \Big( \Big| \sum_{k=1}^n \xi_{k-1} dx_k
\Big|^2 \Big)\Big)} \\ & \le &  4 \lambda  \tau \Big(
\chi_{(\lambda^2/4,\infty)} \Big( 4 \Big| \sum_{k=1}^n \xi_{k-1}
d\alpha_k \Big|^2 \Big)\Big) +  4 \lambda  \tau \Big(
\chi_{(\lambda^2/4,\infty)} \Big( 4 \Big| \sum_{k=1}^n \xi_{k-1}
d\beta_k \Big|^2 \Big)\Big) \\ & + &  4 \lambda  \tau \Big(
\chi_{(\lambda^2/4,\infty)} \Big( 4 \Big| \sum_{k=1}^n \xi_{k-1}
d\gamma_k \Big|^2 \Big)\Big) +  4 \lambda  \tau \Big(
\chi_{(\lambda^2/4,\infty)} \Big( 4 \Big| \sum_{k=1}^n \xi_{k-1}
d\upsilon_k \Big|^2 \Big)\Big) \\ & = & I + II + III +IV.
\end{eqnarray*}
For the first term $I$, we use Chebychev's inequality to deduce:
\begin{equation} \label{I}
I \le \frac{64}{\lambda} \Big\| \sum_{k=1}^n \xi_{k-1} d\alpha_k
\Big\|_2^2 = \frac{64}{\lambda}  \sum_{k=1}^n \|\xi_{k-1}
d\alpha_k \|_2^2 \le \frac{64}{\lambda} \sum_{k=1}^n
\|d\alpha_k\|_2^2 \leq \mathrm{c} \|x\|_1.
\end{equation}
For the second term  $II$, we proceed similarly,
\begin{eqnarray} \label{II}
II & = & 4 \lambda \, \tau \Big( \chi_{(\lambda/2,\infty)} \Big( 2
\Big| \sum_{k=1}^n \xi_{k-1} d\beta_k \Big| \Big)\Big) \\
\nonumber & \le & 16 \Big\| \sum_{k=1}^n \xi_{k-1} d\beta_k
\Big\|_1 \le 16 \sum_{k=1}^n \left\|d\beta_k\right\|_1 \le
\mathrm{c} \|x\|_1.
\end{eqnarray}
For $III$, we note that $|\sum_{k=1}^n \xi_{k-1}d\gamma_k|^2$ is
supported by the projection $$\bigvee_{k\geq 1} \mathrm{supp}
|d\gamma_k|.$$ With this observation, it follows that
\begin{equation} \label{III}
III \leq 4\lambda \T\Big(\bigvee_{k\geq 1} \mathrm{supp}
|d\gamma_k|\Big) \leq \mathrm{c} \|x\|_1.
\end{equation}
For the last term $IV$, we remark  first that since $\xi_{k-1}$
commutes with $d\upsilon_k$, we have $$d\upsilon^*_k
\xi^*_{k-1}=\xi^*_{k-1}d{\upsilon}_k^*.$$ Therefore $$IV = 4
\lambda \tau \Big( \chi_{(\lambda^2/4,\infty)} \Big( 4 \Big|
\sum_{k=1}^n d\upsilon^*_k \xi_{k-1}^* \Big|^2 \Big)\Big)= 4
\lambda  \tau \Big( \chi_{(\lambda^2/4,\infty)} \Big( 4 \Big|
\sum_{k=1}^n \xi_{k-1}^*d\upsilon^*_k  \Big|^2 \Big)\Big). $$
Using the same argument as in $III$, we can conclude that
\begin{equation} \label{IV}
IV \leq \mathrm{c} \|x\|_1.
\end{equation}
Inequality~(\ref{transform-Ineq}) follows immediately from
(\ref{I}, \ref{II}, \ref{III}, \ref{IV}). The proof is complete.
 \fin

\subsection{Non-commutative Burkholder inequality on square
functions}

This subsection is devoted to the non-commutative extension  of
the weak type $(1,1)$ boundedness of square functions of classical
martingales \cite{B2}. In \cite{R3}, the following non-commutative
extension of Burkholder's result was obtained:

\begin{theorem}\label{weak}
There exists an absolute  constant $\mathrm{c} >0$ such that for
any given martingale $x = (x_n)_{n \ge 1}$ that is bounded in
$L_1(\M,\T)\cap L_2(\M,\T)$,  there exist two martingales $y$ and
$z$ with $x = y+z$ and: $$\Big\| \Big(\sum^\infty_{n=1} |dy_n|^2
\Big)^{{1}/{2}} \Big\|_{1,\infty} + \Big\| \Big(\sum^\infty_{n=1}
|dz_n^*|^2 \Big)^{1/2} \Big\|_{1,\infty} \leq \mathrm{c}
\|x\|_1.$$
\end{theorem}

We also refer to \cite{R3} for some applications of Theorem
\ref{weak}. Our purpose is to highlight that most of the
complicated estimates from the proof of Theorem~\ref{weak} in
\cite{R3} can be explained through the decomposition in Section~2.
Namely, we will use the decomposition (\ref{Decomposition2}).

\vskip5pt

\dem First, we recall from \cite{R3} that the general case can be
deduced easily from the special case where $x$ is a positive
martingale. Therefore, without loss of generality, we shall assume
in what follows that $x$ is a positive $L_1$-bounded martingale
and with norm $ \|x\|_1 = 1$. Moreover, we shall only present the
special case where $\M$ is a finite von Neumann algebra with the
trace $\T$ being normalized. The semifinite case only requires
slight changes, see \cite{R3} for further details.

\vskip5pt

We begin by recalling the construction of the martingales $y$ and
$z$ from \cite{R3}. We consider collections of sequences of
pairwise disjoint projections as follows: for $n\geq 1$, set
\begin{equation}
\begin{cases}
p_{0,n} &:=\displaystyle{ \bigwedge^\infty_{k=0}q_{n}^{(2^k)}},
\quad \text{and} \\ p_{i,n} &:=\displaystyle{
\bigwedge^\infty_{k=i}q_{n}^{(2^k)} - \bigwedge^\infty_{k=i-1}
q_{n}^{(2^k)} \quad \text{for $i\geq 1$}},
\end{cases}
\end{equation}
where $(q^{(s)}_n)_{n \ge 0}$ denotes the sequence of Cuculescu's
projections relative to $x$ and $s>0$. The martingales $y$ and $z$
are defined from their respective martingale difference sequences
as follows:
\begin{equation}\label{mainequation}
\begin{cases} dy_1 &:=\displaystyle{ \sum^\infty_{j=0}
\sum_{i \leq j} p_{i,1} dx_1 p_{j,1}};  \\ dy_k &:=\displaystyle{
\sum^\infty_{j=0}\sum_{i \leq j} p_{i,k-1} dx_k p_{j,k-1}}  \quad
\text{for $k\geq 2$}; \\  dz_1 &:=
\displaystyle{\sum^\infty_{j=0}\sum_{i > j} p_{i,1} dx_1 p_{j,1}};
\\ dz_k &:= \displaystyle{\sum^\infty_{j=0}\sum_{i > j}
p_{i,k-1} dx_k p_{j,k-1}} \quad \text{for $k\geq 2$}.
\end{cases}
\end{equation}
We refer to \cite{R3} for the fact that indeed $(dy_k)_{k \ge 1}$
and $(dz_k)_{k \ge 1}$ are martingale difference sequences.
Moreover, as already explained in \cite{R3}, it suffices to see
that there exists an absolute constant $\mathrm{c}$ such that for
every non-negative integer $m$,
\begin{equation}\label{maininequality}
 2^m \, \tau \Big( \chi_{(2^m,
\infty)} \big( \mathcal{S}_C(y) \big) \Big) \le \mathrm{c}.
\end{equation}
This will be done in two steps:

\vskip5pt

\noindent \textbf{Step 1.} In what follows, we take $\lambda =
2^m$ for some $m \ge 0$. We begin as in \cite{R3} by truncating
the series which defines $dy_k$. More concretely, applying Lemma
\ref{Quasi-Triangle} and Proposition \ref{Cuculescu} we make the
following reduction, see \cite{R3} for further details.

\begin{proposition} \label{PropositionA}
\rm{\cite[Proposition~A]{R3}.} We have
\begin{equation*}
\lambda \T \Big( \chi_{(\lambda,\infty)}(\mathcal{S}_{C}(y)) \Big)
\leq 2 \lambda \T
\Big(\chi_{(\lambda/2,\infty)}(\mathcal{S}_{C}^{(m)}(y)) \Big) +
4,
\end{equation*}
where $\mathcal{S}_{C}^{(m)}(y)$ denotes the following square
function $$\mathcal{S}_{C}^{(m)}(y) = \left( \Big| \sum_{j=0}^{m}
\sum_{i\leq j} p_{i,1} dx_1 p_{j,1} \Big|^2 + \sum_{k=2}^{\infty}
\Big| \sum_{j=0}^{m} \sum_{i\leq j} p_{i,k-1} dx_k p_{j,k-1}
\Big|^2 \right)^{1/2}.$$
\end{proposition}

\noindent \textbf{Step 2.} According to (\ref{maininequality}) and
Proposition~\ref{PropositionA}, we need to estimate $$2 \lambda \T
\Big(\chi_{(\lambda/2,\infty)}(\mathcal{S}_{C}^{(m)}(y)) \Big).$$
For this, we consider the  decomposition (\ref{Decomposition2}) of
$x$ associated to the parameter $\lambda = 2^m$. This gives $x =
\alpha' + \beta' + \gamma' + \upsilon'$ as in Corollary
\ref{Gundy}. For $k \geq 1$, set
\[
\mathsf{P}^{(m)}_k :=(p_{i,k})_{i=0}^m.
\]
Then with this notation,
\[
\mathcal{S}_{C}^{(m)}(y) = \Big( \big|
\cal{T}^{\mathsf{P}^{(m)}_1}(dx_1) \big|^2 + \sum_{k=2}^\infty
\big| \cal{T}^{\mathsf{P}^{(m)}_{k-1}}(dx_k) \big|^2 \Big)^{1/2}.
\]
We make the following crucial observation:

\begin{lemma}
The martingales $\gamma'$ and $\upsilon'$ do not contribute to the
quantity $\mathcal{S}_{C}^{(m)}(y)$. In particular,
\[\mathcal{S}_{C}^{(m)}(y) = \left( \big| \cal{T}^{\mathsf{P}^{(m)}_1}
(d\alpha'_1 + d\beta'_1) \big|^2 + \sum_{k=2}^\infty
|\cal{T}^{\mathsf{P}^{(m)}_{k-1}}(d\alpha'_k +d\beta'_k)) \big|^2
\right)^{1/2}.\]
\end{lemma}

\dem Recall from (\ref{Decomposition2}) that $d\gamma'_1
+d\upsilon'_1=0$ and for every $k\geq 2$,
\[
d\gamma'_k  + d\upsilon'_k= dx_k - q^{(2^m)}_{k-1} dx_k
q^{(2^m)}_{k-1}.
\]
The key point here is that
\[
p_{i,k-1}( d \gamma'_k + d\upsilon'_k) p_{j,k-1} = 0 \quad
\mbox{for any} \quad 0 \le i,j \le m \quad \mbox{whenever} \quad k
\ge 2.
\]
This gives $\cal{T}^{\mathsf{P}^{(m)}_{k-1}}(d\gamma'_k
+d\upsilon'_k)=0$ for every $k\geq 2$, which proves the lemma.
\fin

Now, using the elementary inequality for operators $|a+b|^2 \le 2
|a|^2 + 2 |b|^2$ and Lemma~\ref{Quasi-Triangle} above, we deduce
that
\begin{eqnarray*}
\lefteqn{2\lambda \T \Big(\chi_{(\lambda/2,
\infty)}(\mathcal{S}_{C}^{(m)}(y)) \Big)} \\ & \leq & 4\lambda \T
\left( \chi_{(\lambda^2/8, \infty )} \Big( 2 \big|
\cal{T}^{\mathsf{P}^{(m)}_1}(d\alpha'_1) \big|^2 + 2
\sum_{k=2}^\infty \big|
\cal{T}^{\mathsf{P}^{(m)}_{k-1}}(d\alpha'_k) \big|^2 \Big) \right)
\\ & + &
4\lambda \T \left( \chi_{(\lambda^2/8, \infty)} \Big( 2 \big|
\cal{T}^{\mathsf{P}^{(m)}_1}(d\beta'_1) \big|^2 + 2
\sum_{k=2}^\infty \big|
\cal{T}^{\mathsf{P}^{(m)}_{k-1}}(d\beta'_k) \big|^2 \Big) \right)
\\ & = & I + II.
\end{eqnarray*}
Chebychev's inequality gives
\begin{eqnarray*}
I & \le & \frac{64}{\lambda} \Big( \big\|
\cal{T}^{\mathsf{P}^{(m)}_1}(d\alpha'_1) \big\|_2^2 +
\sum_{k=2}^{\infty} \big\|
\cal{T}^{\mathsf{P}^{(m)}_{k-1}}(d\alpha'_k) \big\|_2^2 \Big).
\end{eqnarray*}
Since triangular truncations are  orthogonal projections in
$L_2(\mathcal{M},\T)$, we deduce
\[
I \le \frac{64}{\lambda} \sum_{k=1}^\infty\|d\alpha'_k\|_2^2  \leq
\mathrm{c}.
\]
where the last inequality follows from Corollary~\ref{Gundy} (i).
For $II$, we have
\[
II \leq 16 \Big\| \Big( \big|
\cal{T}^{\mathsf{P}^{(m)}_1}(d\beta'_1) \big|^2 +
\sum_{k=2}^{\infty} \big|
\cal{T}^{\mathsf{P}^{(m)}_{k-1}}(d\beta'_k) \big|^2 \Big)^{1/2}
\Big\|_{1,\infty}.
\]
Therefore we can apply Lemma \ref{Truncation} to obtain that
\[
II \le \mathrm{c}  \sum_{k=1}^{\infty} \|d\beta'_k\|_1 \le
\mathrm{c}.
\]
The last estimate follows once more  from Corollary~\ref{Gundy}
(ii). Thus, combining the estimates for $I$ and $II$ with
Proposition~\ref{PropositionA}, the desired inequality
(\ref{maininequality}) follows. \fin

\begin{remark}
\emph{Let us consider a positive non-commutative martingale $x$
and let $x = \alpha' + \beta' + \gamma' + \upsilon'$ stand for the
decomposition given in Corollary \ref{Gundy} associated to
$\lambda$. Then, using the properties stated in Corollary
\ref{Gundy}, it is not difficult to check that}
\begin{eqnarray*}
\max \Big\{ \lambda \tau \big( \chi_{(\lambda,\infty)}
(\mathcal{S}_R(\alpha')) \big), \, \lambda \tau \big(
\chi_{(\lambda,\infty)} (\mathcal{S}_C(\alpha')) \big) \Big\} &
\le & \mathrm{c} \|x\|_1, \\ \max \Big\{ \lambda \tau \big(
\chi_{(\lambda,\infty)} (\mathcal{S}_R(\beta')) \big), \, \lambda
\tau \big( \chi_{(\lambda,\infty)} (\mathcal{S}_C(\beta')) \big)
\Big\} & \le & \mathrm{c} \|x\|_1, \\ \max \Big\{ \lambda \tau
\big( \chi_{(\lambda,\infty)} (\mathcal{S}_R(\upsilon')) \big), \,
\lambda \tau \big( \chi_{(\lambda,\infty)}
(\mathcal{S}_C(\gamma')) \big) \Big\} & \le & \mathrm{c} \|x\|_1.
\end{eqnarray*}
\emph{Unfortunately, in contrast with the commutative case,
Theorem \ref{Main} does not follow automatically from these
estimates. Namely, in the commutative case for any given
martingale $x$ and any $\lambda > 0$, we can take Gundy's
decomposition associated to $\lambda$ and then the commutative
analogue of the estimates given above (where row and column square
functions coincide) provides the desired weak type inequality, see
\cite{G} for the details. However, in the non-commutative setting
we are forced to decompose the martingale $x$ into two other
martingales $x=y+z$ \emph{before} being able to apply Gundy's
decomposition. This is justified by the fact that the
non-commutative weak Hardy space
$\mathcal{H}_{1,\infty}(\mathcal{M}, \tau)$ is the sum of two
quasi-Banach spaces, see \cite{PX1,R3} for further details. Thus,
the classical proof of Theorem \ref{Main} for commutative
martingales does not work here since we would have to consider the
Gundy's decomposition of $y$ and $z$ separately. However, even in
the case where both $y$ and $z$ were $L_1$-bounded, this only
allows us to control our terms by the norms of $y$ and $z$ in
$L_1(\mathcal{M})$, not by the norm of $x$ in $L_1(\mathcal{M})$.}
\end{remark}

\subsection{Co-lacunary sequences in non-commutative $L_1$-spaces}

Let $\mathrm{X}$ be a Banach space. A sequence $(x_n)_{n \ge 1}$
in $\mathrm{X}$ is said to be 2\emph{-co-lacunary} if there is
$\delta>0$ such that for any finite sequence $(a_n)_{n \ge 1}$ of
scalars,
\[
\delta \Big( \sum_{n \ge 1} |a_n|^2 \Big)^{1/2} \le \Big\| \sum_{n
\ge 1} a_n x_n \Big\|_{\mathrm{X}}.
\]
This property can also be described by saying that $(x_n)_{n \ge
1}$ dominates the unit vector basis of $l_2$, but we will follow
the term 2-co-lacunary from \cite{ALFR} which was motivated by the
terminology {\it lacunary sequences} used in \cite{KP} for a dual
pro\-perty. In \cite{ALFR}, Aldous and Fremlin proved the
remarkable result that if $(\Omega,\cal{F}, \mu)$ is a
pro\-bability space, then every uniformly integrable martingale
difference sequence which is bounded away from zero is
$2$-co-lacunary in $L_1(\Omega,\mu)$. Using such result, they
deduced the following subsequence principle in $L_1$-spaces: every
bounded sequence in $L_1(\Omega,\mu)$ has either a convergent or a
$2$-co-lacunary subsequence in $L_1(\Omega,\mu)$.

\vskip5pt

The principal result of this section is Theorem~\ref{main-colac}
below which extends the main result of Aldous and Fremlin in
\cite{ALFR} on classical martingale difference sequences to the
non-commutative setting. The proof in \cite{ALFR} is very involved
and based on several use of stopping times. Another proof was also
given by Dor in \cite{Dor}. Our proof below uses the decomposition
(\ref{Decomposition2}). This approach seems to be overlooked for
the commutative case.

\begin{theorem}\label{main-colac}
Let $(d_k)_{k \ge 1}$ be a martingale difference sequence in
$L_1(\M,\T)$  with:
\begin{itemize}
\item[(i)] $\gamma = \inf \big\{\|d_k\|_1 \, \big| \ k \geq 1
\big\} >0$;
\item[(ii)] $\big\{ d_k \, \big| \ k \geq 1 \big\}$ is relatively
weakly compact
in $L_1(\M,\T)$.
\end{itemize}
Then, the sequence $(d_k)_{k \ge 1}$ is a $2$-co-lacunary sequence
in $L_1(\M,\T)$.
\end{theorem}

From the closed graph theorem, Theorem~\ref{main-colac} can be
reformulated as follows:

\addtocounter{theorem}{-1}
\renewcommand{\thetheorem}{\arabic{section}.\arabic{theorem}'}

\begin{theorem}
Let $(d_k)_{k \ge 1}$ be the  sequence described in Theorem
$\ref{main-colac}$. Suppose that $(a_k)_{k \ge 1}$ is a sequence
of scalars such that the series $\sum_{k \ge 1} a_kd_k$ is
convergent in $L_1(\M,\T)$. Then we have $$\sum_{k \ge 1} |a_k|^2
<\infty.$$
\end{theorem}

\renewcommand{\thetheorem}{\arabic{section}.\arabic{theorem}}

\dem The proof will be divided into several cases.

\vskip5pt

\noindent{\bf Case~A:} Assume first that $(d_k)_{k \ge 1}$ is a
sequence of self-adjoint operators and $(a_k)_{k \ge 1}$ is a
sequence in $\mathbb{R}$. We start by noting that since $(d_k)_{k
\ge 1}$ is relatively weakly compact (and thus, equiintegrable in
the sense of \cite{Ran10}), condition~(i) of
Theorem~\ref{main-colac} is equivalent to:
\begin{equation}\label{sigma}
\sigma :=\inf\left\{\left\|d_k\right\|_{L_1(\M,\T) + \M} \, \big|
\ \ k \geq 1\right\} >0.
\end{equation}
Indeed, if $\liminf_{k\to \infty}\|d_k\|_{L_1(\M,\T) + \M}=0$,
then there is a subsequence $(d_{k_j})_{j \ge 1}$ that converges
to zero in $L_1(\M,\T) +\M$, and a fortiori, it converges to zero
in the measure topology. By \cite[Proposition~2.11]{Ran10}, we
have $\lim_{j \to \infty}\|d_{k_j}\|_1 = 0$ which violates
condition~(i).

\vskip5pt

Next, we observe that $(a_kd_k)_{k \ge 1}$ is a (self-adjoint)
martingale difference sequence and denote by $y=(y_n)_{n \ge 1}$
the corresponding martingale. By assumption, $y$ is a
$L_1$-bounded self-adjoint martingale. For every $\lambda>0$, we
can consider the sequence of Cuculescu's projections associated to
the martingale $y$ and $\lambda > 0$. We claim that
\[
\lim_{\lambda \to \infty} \sup \Big\{ \big\|d_k \big( {\bf 1} -
q_k^{(\lambda)} \big) \big\|_{L_1(\M,\T) +\M} \, \big| \ k \in
\mathbb{N} \Big\} = 0.
\]
Indeed, for every $k \in \mathbb{N}$ and $\lambda>0$, we have (see
e.g. \cite[Proposition~2.a.2]{LT}):
\begin{equation*}
\big\|d_k \big( {\bf 1}-q_k^{(\lambda)} \big) \big\|_{L_1(\M,\T)
+\M} = \int_0^1 \mu_t \big( d_k \big( {\bf 1}-q_k^{(\lambda)}
\big) \big) \ dt.
\end{equation*}
Using properties of singular value functions from \cite{FK},
\begin{eqnarray*}
\big\|d_k \big( {\bf 1}-q_k^{(\lambda)} \big) \big\|_{L_1(\M,\T)
+\M} & \leq & \int_0^1 \mu_t(d_k) \, \mu_t \big( {\bf 1} -
q_k^{(\lambda)} \big) \ dt
\\ & \leq & \int_0^1 \mu_t(d_k) \, \chi_{[0,\T({\bf
1}-q_k^{(\lambda)})]}(t) \ dt \\ & \leq & \int_0^{1 \wedge
\lambda^{-1}\|y\|_1} \mu_t(d_k) \ dt.
\end{eqnarray*}
Since $\{d_k \, | \ k \geq 1\}$ is relatively weakly compact in
$L_1(\M,\T)$, it is a fortiori relatively $\sigma(L_1(\M, \T) +\M,
L_1(\M, \T) \cap \M)$-compact in $L_1(\M,\T) +\M$. According to
the terminology used in Section \ref{Section1}, we note that
$L_1(\M,\T) \cap \M \subset S_0(\M,\T)$. Then, by
Theorem~\ref{weakcom2} $\{\mu(d_k) \, | \ k \geq 1\}$ is
relatively $\sigma(L_1 +L_\infty, L_1 \cap L_\infty)$-compact in
$L_1[0,\T({\bf 1}))+L_\infty[0,\T({\bf 1}))$ and therefore
$\{\mu(d_k)\chi_{[0,1]} \, | \ k \geq 1\}$ is relatively weakly
compact in $L_1[0,1]$. In other words, it is uniformly integrable
in $L_1[0,1]$. Thus we conclude that $$\lim_{\lambda \to \infty}
\sup \left\{ \int_0^{1 \wedge \lambda^{-1}\|y\|_1} \mu_t(d_k) \ dt
\ \big| \ k \ge 1 \right\} = 0,$$ which proves the claim. For the
remainder of the proof, we fix $\lambda> 0$ so that,
\begin{equation}\label{eta}
\sup \Big\{ \big\|d_k({\bf 1} - q_k^{(\lambda)})
\big\|_{L^1(\M,\T) +\M} \, \big| \ k \ge 1 \Big\} \leq \sigma/5
\end{equation}
where $\sigma$ is from (\ref{sigma}). We will simply write
$(q_n)_{n \ge 0}$ for $(q_n^{(\lambda)})_{n \ge 0}$. We note from
(\ref{eta}) that for every    $k\geq 1$,
\begin{eqnarray*}
\sigma & \leq & \|d_k\|_{L_1(\M,\T)+\M} \\ & \leq & \big\|q_k d_k
q_k - \mathsf{E}_{k-1}(q_k d_k q_k) \big\|_{L_1(\M,\T)+\M} \\ & +
& \big\|q_k d_k ({\bf 1}-q_k) - \mathsf{E}_{k-1}(q_k d_k ({\bf
1}-q_k )) \big\|_{L_1(\M,\T)+\M} \\ & + & \big\|({\bf 1}-q_k)d_k -
\mathsf{E}_{k-1}(({\bf 1}-q_k)d_k ) \big\|_{L_1(\M,\T)+\M} \\ &
\leq & \big\| q_k d_k q_k - \mathsf{E}_{k-1}(q_k d_k q_k) \big\|_2
\\ & + & 2 \big\|q_k d_k ({\bf 1}-q_k)\big\|_{L_1(\M,\T) +\M} \\ & + & 2
\big\| ({\bf 1}-q_k)d_k \big\|_{L_1(\M,\T) +\M} \\ & \leq & \big\|
q_k d_k q_k - \mathsf{E}_{k-1}(q_k d_k q_k) \big\|_2 + 4 \sigma/5.
\end{eqnarray*}
Therefore,
\begin{equation}\label{alpha2}
\inf \Big\{ \big\|q_k d_k q_k - \mathsf{E}_{k-1}(q_k d_k q_k)
\big\|_2 \, \big| \ \ k \geq 1 \Big\} \geq \sigma/5.
\end{equation}
Let us consider the decomposition $y=\alpha' +\beta' +\gamma'
+\upsilon'$ of the martingale $y$ according to
Corollary~\ref{Gundy} and Remark \ref{Non-positive} and relative
to the parameter $\lambda > 0$ fixed above. Then we have
$\|\alpha'\|_2^2 \leq \mathrm{c} \lambda\|y\|_1$. Recall that for
every $k\geq 1$, $$d\alpha'_k =q_kdy_kq_k
-\mathsf{E}_{k-1}(q_kdy_k q_k) =a_k\left(q_kd_kq_k
-\mathsf{E}_{k-1}(q_kd_k q_k)\right).$$ This gives,
\begin{equation*}
\sum_{k\geq 1} \left|a_k\right|^2 \big\|q_kd_kq_k
-\mathsf{E}_{k-1}(q_kd_k q_k)\big\|_2^2 \leq \mathrm{c}
\lambda\left\|y\right\|_1,
\end{equation*}
and therefore by (\ref{alpha2}) we conclude,
\begin{equation*}
\sigma^2 \sum_{k\geq 1} \left|a_k\right|^2 \leq 25 \mathrm{c}
\lambda\left\|y\right\|_1 <\infty.
\end{equation*}
The proof for this case is complete.

\vskip5pt

\noindent{\bf Case~B:} Assume now that $(d_k)_{k \ge 1}$ is not
necessarily a sequence of self-adjoint operators and $(a_k)_{k \ge
1}$ is a sequence in $\mathbb{R}$. Consider the (semifinite) von
Neumann algebra $\M \oplus_\infty \M$ with the trace
$\widetilde{\T}=\T \oplus_\infty \T$ and  the filtration $(\M_n
\oplus_\infty \M_n)_{n \ge 1}$. For $k \geq 1$, let
\[
\widetilde{d}_k := \Big( (d_k +d_k^*)/2, (d_k-d_k^*)/2i \Big) \in
\M_k \oplus_\infty \M_k.
\]
Then $(\widetilde{d}_k)_{k \ge 1}$ is a self-adjoint martingale
difference sequence in $L_1(\M \oplus_\infty \M,\widetilde{\T})$
that clearly verifies the conditions (i) and (ii) of
Theorem~\ref{main-colac}. If $\sum_{k \ge 1} a_kd_k$ is
convergent, then so is the series $\sum_{k \ge 1} a_kd_k^*$.
Consequently, the series $$\sum_{k \ge 1} a_k\widetilde{d}_k$$
converges in $L_1(\M \oplus_\infty \M,\widetilde{\T})$. Hence,
from Case~A, we get $\sum_{k \ge 1}|a_k|^2 <\infty$.

\vskip5pt

\noindent{\bf Case~C:} For the general case where $a_k \in \C$, we
set
\begin{equation*}
\gamma_k :=
\begin{cases}
a_k/|a_k| \ \ & \text{if} \ a_k \neq 0 \\ 1  \ \ \ &\text{if} \
a_k=0.
\end{cases}
\end{equation*}
Then it is clear that $(\hat{d}_k)_{k \ge 1} =(\gamma_k d_k)_{k
\ge 1}$ is a martingale difference sequence that satisfies
conditions (i) and (ii) of Theorem~\ref{main-colac}. Moreover,
since we assume by hypothesis that the series $$\sum_{k \ge 1} a_k
d_k = \sum_{k \ge 1} |a_k|\hat{d}_k$$ is convergent, Case~B
insures that $\sum_{k \ge 1} |a_k|^2 <\infty$. The proof is
complete. \fin

As a consequence of Theorem~\ref{main-colac}, we have the
following result that generalizes a result from \cite{ALFR} (see
also \cite{Dor} for a quantitative version) to non-commutative
spaces.

\begin{corollary}\label{fhyper}
Assume that $\M$ is semifinite and hyperfinite. Let $(x_n)_{n \ge
1}$ be a bounded sequence in $L_1(\M,\T)$. Then either $(x_n)_{n
\ge 1}$ has a convergent subsequence or it has a $2$-co-lacunary
subsequence.
\end{corollary}

For the proof, we will use the following perturbation lemma from
\cite{ALFR}.

\begin{lemma}\label{perturbation}
Let $\mathrm{X}$ be a normed space and $(x_n)_{n \ge 1}$ be a
bounded sequence in $\mathrm{X}$. Then the following properties
hold:
\begin{itemize}
\item[(a)] If $(x_n)_{n \ge 1}$ is $2$-co-lacunary and $x \in
\mathrm{X}$, then there exists $m\in \mathbb{N}$ such that $(x_n
-x )_{n \ge m}$ is $2$-co-lacunary;
\item[(b)] If $(x_n)_{n \ge 1}$ is $2$-co-lacunary and
$(y_n)_{n \ge 1} \subset \mathrm{X}$ with $\sum_{n \ge 1}
\|x_n-y_n\|_{\mathrm{X}}$ being convergent, then there exists
$m\in \mathbb{N}$ such that $(y_n)_{n \ge m}$ is $2$-co-lacunary.
\end{itemize}
\end{lemma}

\noindent{\bf Proof of Corollary~\ref{fhyper}.} Assume that
$(x_n)_{n \ge 1}$ has no convergent subsequence. By Rosenthal's
$l_1$-theorem, either $(x_n)_{n \ge 1}$ has a subsequence
equivalent to the unit vector basis of $l_1$, and therefore
$2$-co-lacunary (in fact, $1$-co-lacunary), or $(x_n)_{n \ge 1}$
has a  weakly convergent subsequence. Assume w.l.o.g. that
$(x_n)_{n \ge 1}$ converges to $x$ weakly. Then, according to
Lemma \ref{perturbation} (a), it suffices to show that $(x_n-x)_{n
\ge 1}$ has a $2$-co-lacunary subsequence.

\vskip5pt

If $\cal M$ is hyperfinite then $\cal M = \overline{\cup_\alpha
\cal M_\alpha}$ (weak* closure) where $(\M_\alpha)_{\alpha \in I}$
is a net of finite dimensional *-subalgebras directed by
inclusion. There exist contractive projections $E_\alpha: \cal M
\to \M_\alpha$ which are simultaneously contractions from $\M \to
\M_\alpha$ and $L_1 (\cal M, \tau) \to L_1 (\cal
M_{\alpha},\tau_{\alpha})$, where $\tau_{\alpha}$ denotes the
restriction of $\tau$ on $\cal M_{\alpha}$. The projections
$E_\alpha$'s satisfy $E_\alpha=E_\alpha E_\beta$ for $\alpha \leq
\beta$. Moreover, for every $f \in L_1(\cal M, \tau)$,
$\lim_{\alpha} \Vert E_{\alpha} (f)-f\Vert_1=0$. For $n\geq 1$,
let $f_n=x_n-x$. Then by assumption $(f_n)_{n \ge 1}$ is a weakly
null sequence in $L_1(\M,\T)$. Moreover, since $(f_n)_{n \ge 1}$
does not converges in the norm of $L_1(\mathcal{M}, \tau)$, we may
assume w.l.o.g. that the sequence $(f_n)_{n \ge 1}$ itself
satisfies
\begin{equation} \label{inf>0}
\inf \Big\{ \|f_n\|_1 \, \big| \ n \ge 1 \Big\}
> 0.
\end{equation}
Set $n_1 = 1$, choose $\alpha_1 \in I$ such that
\begin{equation} \label{alpha1}
\big\| f_{1} - E_{\alpha_{1}} (f_{1}) \big\|_1 < {2}^{-2}.
\end{equation}
Since $(f_{n})_{n \ge 1}$ converges weakly to zero and
$\M_{\alpha_1}$ is finite dimensional, we have
\[
\lim_{n \to \infty} \left\Vert E_{\alpha_{1}} (f_{n})\right\Vert_1
= 0.
\]
Choose $n_2 > n_1 = 1$ such that,
\[
\left\|E_{\alpha_{1}} (f_{n_2}) \right\|_1 <{2^{-3}}
\]
and $\alpha_2 > \alpha_1$ so that,
\[ \big\| f_{n_{2}} -
E_{\alpha_{2}} (f_{n_{2}}) \big\|_1< {2^{-3}}.
\]
Inductively, one gets a sequence $(n_{k})_{k \ge 1} \subseteq
\mathbb {N}$ and $\alpha_1 < \alpha_2 < \cdots < \alpha_k <
\cdots$ in $I$ such that for every $k \geq 2$,
\begin{equation} \label{alphak}
\max \Big\{ \left\|E_{\alpha_{k-1}} (f_{n_{k}}) \right\|_1,  \,
\big\| f_{n_{k}} - E_{\alpha_{k}} (f_{n_{k}}) \big\|_1 \Big\} <
{2^{-(k+1)}}.
\end{equation}
For $k \geq 2$, set
\begin{equation*}\begin{cases}
v_1 &:= E_{\alpha_1}(f_1), \cr v_k &:= E_{\alpha_{k}}(f_{n_{k}}) -
E_{\alpha_{k-1}} (f_{n_{k}}).
\end{cases}\end{equation*}

\begin{lemma}\label{mds2}
The sequence $(v_k)_{k \ge 1}$ satisfies:
\begin{itemize}
\item[(a)] $\big\| f_{n_k}-v_k \big\|_1\leq 2^{-k}$ for every
$k \geq 1$;
\item[(b)] $\liminf_{k \to \infty} \|v_k\|_1 > 0$;
\item[(c)] $\big\{ v_k \, | \ k\geq 1 \big\}$ is relatively weakly
compact in $L_1(\M,\T)$.
\end{itemize}
\end{lemma}

\dem Property (a) follows from (\ref{alpha1}, \ref{alphak}) and
the triangle inequality. Property (b) follows from (a) and
(\ref{inf>0}). Finally, (c) follows directly from (a) and the fact
that $(f_{n_k})_{k \ge 1}$ is weakly null. \fin

 \vskip5pt

Note that if  $\M$ is finite, then  the projections $E_\alpha$'s
defined above are conditional expectations so the sequence
$(v_k)_{k \ge 1}$ is clearly a martingale difference sequence with
respect to the filtration $(\M_{\alpha_k})_{k \ge 1}$. Assume that
$\M$ is infinite. For every $k \geq 1$, let $p_k$ be the
self-adjoint projection which is the unit element of
$\M_{\alpha_{k}}$. When $\tau$ is infinite, $p_k \neq {\bf 1}$.
The contractive projection $E_{\alpha_{k}}$ can be written as
$E_{\alpha_{k}} =\cal{E}_k \circ Q_{\alpha_{k}}$ where
$Q_{\alpha_{k}}(x)=p_kxp_k$ and $\cal{E}_k$ is the unit preserving
conditional expectation from  the finite von Neumann algebra
$p_k\M p_k$ onto  the von Neumann subalgebra $\M_{\alpha_{k}}$.
Since $\big\{ v_k \, | \ k\geq 1 \big\}$ is relatively weakly
compact in $L_1(\M,\T)$ and $(p_{k}-p_{k-1})_{k\geq 2}$ is a
disjoint sequence of projections, we have (see \cite{Ran10})
\[
\lim_{k \to \infty}
\big\|(p_{k}-p_{k-1})v_k(p_{k}-p_{k-1})\big\|_1=0.
\]
By taking subsequence if necessary, we can assume that for every
$k\geq 2$,
\begin{equation}\label{v}
\big\|(p_{k}-p_{k-1})v_k(p_{k}-p_{k-1})\big\|_1 \leq 2^{-k}.
\end{equation}
For $k \geq 2$, set
\begin{equation*}\begin{cases}
w_1 &:= v_1, \cr w_k &:=v_k -(p_{k}-p_{k-1})v_k(p_{k}-p_{k-1}).
\end{cases}\end{equation*}
\begin{lemma}\label{mds}
The sequence $(w_k)_{k \ge 1}$ is a martingale difference sequence
with:
\begin{itemize}
\item[(a)] $\big\| f_{n_k}-w_k \big\|_1\leq 2^{-(k-1)}$ for every
$k \geq 1$;
\item[(b)] $\liminf_{k \to \infty} \|w_k\|_1 > 0$;
\item[(c)] $\big\{ w_k \, | \ k\geq 1 \big\}$ is relatively weakly
compact in $L_1(\M,\T)$.
\end{itemize}
\end{lemma}
\dem Properties (a), (b), and (c) follow directly from
Lemma~\ref{mds2} and (\ref{v}). Hence it remains to show that
$(w_k)_{k\geq 1}$ is a martingale difference sequence. For $k\geq
1$, we set
\[
\cal{S}_k:=\M_{\alpha_{k}} + \, \sum_{s \ge k}
(p_{s+1}-p_s)\M_{\alpha_{s+1}}(p_{s+1}- p_s).
\]
Clearly, $(\cal{S}_k)_{k \ge 1}$ is an increasing sequence of von
Neumann subalgebras of $\M$ and, by the definition of the
$\mathcal{M}_{\alpha}$'s, we may assume that $\M=\overline{\cup_{k
\ge 1} \cal{S}_k}$ (weak* closure). Define $\mathbb{E}_k: \M \to
\cal{S}_k$ by setting
\[
\mathbb{E}_k(x) = E_{\alpha_{k}}(x) + \sum_{s\geq k} (p_{s+1}-
p_s)E_{\alpha_{s+1}}(x)(p_{s+1} -p_s),
\]
where the sum is taken with respect to the weak* topology. The
right hand side is well defined since the series $\sum_{s\geq k}
(p_{s+1}- p_s)E_{\alpha_{s+1}}(x)(p_{s+1} -p_s)$ is weakly
unconditionally Cauchy in $\M$. In fact, for any finite set $S
\subset [k, \infty)$, $$\Big\| \sum_{s\in S} (p_{s+1}- p_s)
E_{\alpha_{s+1}}(x)(p_{s+1} -p_s) \Big\|_\infty \leq \left\Vert x
\right\Vert_\infty.$$ The operator $\mathbb{E}_k$ is clearly a
(contractive) positive projection and $\mathbb{E}_k({\bf 1})={\bf
1}$. Thus, according to \cite{Tomi}, the mapping $\mathbb{E}_k$ is
a conditional expectation. Moreover, for every $x \in \M$, we have
\begin{eqnarray*}
\tau(\mathbb{E}_k(x)) & = & \tau(E_{\alpha_{k}}(x)) + \sum_{s\geq
k} \tau \big( (p_{s+1}- p_s)E_{\alpha_{s+1}}(x)(p_{s+1}- p_s)
\big) \\ & = & \tau \big( \cal{E}_k(p_kxp_k) \big) + \sum_{s\geq
k} \tau \big( \cal{E}_{s+1} ((p_{s+1}- p_s)x(p_{s+1} -p_s)) \big)
\\ & = & \tau(p_kxp_k) + \sum_{s\geq k} \tau \big( (p_{s+1}-
p_s)x(p_{s+1} -p_s) \big) = \tau(x).
\end{eqnarray*}
It is clear that for every $k\geq 2$, $E_{\alpha_{k-1}}(v_k)=0$
and $E_{\alpha_{k-1}}((p_k-p_{k-1})v_k(p_k-p_{k-1}))=0$. We can
deduce from   the definition of  $\mathbb{E}_{k-1}$ that
$\mathbb{E}_{k-1}(w_k)=0$. Therefore, $ (w_k)_{k \ge 1}$ is a
martingale difference sequence with respect to the filtration
$(\cal{S}_k)_{k \ge 1}$. The lemma is proved. \fin

From Theorem \ref{main-colac} and Lemma~\ref{mds}, we know that
$(w_k)_{k \ge m_1}$ is $2$-co-lacunary for some $m_1 \ge 1$.
Moreover, since $\sum_{k \ge 1} \|f_{n_k}-w_k\|_1<\infty$, it
follows from Lemma~\ref{perturbation} that there exists $m_2 \geq
m_1$ such that $(f_{n_k})_{k \ge m_2} = (x_{n_k} - x)_{k \ge m_2}$
is $2$-co-lacunary. The proof of Corollary~\ref{fhyper} is
complete. \qed

\begin{remark}
\emph{We do not know if Corollary~\ref{fhyper} is valid without
the hyperfinite assumption. We leave this as a problem for the
interested reader.}
\end{remark}


\bibliographystyle{amsplain}

\begin{thebibliography}{10}

\bibitem {ALFR} D.J. Aldous and D.H. Fremlin, \emph{Colacunary
sequences in {$L$}-spaces}, Studia Math. \textbf{71} (1981/82),
297-304.

\bibitem {B2} D.L. Burkholder, \emph{Martingale transforms}, Ann.
Math. Statist. \textbf{37} (1966), 1494-1504.

\bibitem {B} D.L. Burkholder, \emph{Distribution function
inequalities for martingales}, Ann. Probab. \textbf{1} (1973),
19-42.

\bibitem {BG} D.L. Burkholder and R.F. Gundy, \emph{Extrapolation
and interpolation of quasi-linear operators on martingales}, Acta
Math. \textbf{124} (1970), 249-304.

\bibitem {CZ} A.P. Calder\'{o}n and A. Zygmund, \emph{On the existence
of certain singular integrals}, Acta Math. \textbf{88} (1952),
85-139.

\bibitem {C} I. Cuculescu, \emph{Martingales on von Neumann
algebras}, J. Multivariate Anal. \textbf{1} (1971), 17-27.

\bibitem {DDdP} P.G. Dodds, T.K. Dodds and B. de Pagter,
\emph{Noncommutative Banach function spaces}, Math. Z.
\textbf{201} (1989), 583-597.

\bibitem {DDP3} P.G. Dodds, T.K. Dodds and B. de Pagter,
\emph{Noncommutative {K}\"othe duality}, Trans. Amer. Math. Soc.
\textbf{339} (1993), 717-750.

\bibitem {DSS} P.G. Dodds, F.A. Sukochev, and G.Schl{\"u}chtermann,
\emph{Weak compactness criteria in symmetric spaces of measurable
operators}, Math. Proc. Cambridge Philos. Soc. \textbf{131}
(2001), 363-384.

\bibitem {Dor} L. Dor, \emph{On co-lacunary sequences in
{$L\sb{1}$}}, Bull. London Math. Soc. \textbf{14} (1982), 410-414.

\bibitem {FK} T. Fack and H. Kosaki, \emph{Generalized $s$-numbers
of $\tau$-measurable operators}, Pacific J. Math. \textbf{123}
(1986), 269-300.

\bibitem {GA} A.M. Garsia, Martingale inequalities: {S}eminar
notes on recent progress. W. A. Benjamin, Inc., Reading,
Mass.-London-Amsterdam, 1973, Mathematics Lecture Notes Series.

\bibitem {G} R.F. Gundy, \emph{A decomposition for $L^1$-bounded
martingales}, Ann. Math. Statist. \textbf{39} (1968), 134-138.

\bibitem {G2} R.F. Gundy, \emph{On the class $L \log L$,
martingales, and singular integrals}, Studia Math. \textbf{33}
(1969), 109-118.

\bibitem {J1} M. Junge, \emph{Doob's inequality for non-commutative
martingales}, J. Reine Angew. Math. \textbf{549} (2002), 149-190.

\bibitem {JM} M. Junge and M. Musat, \emph{A non-commutative
version of the John-Nirenberg theorem}. Preprint 2004.

\bibitem {JX} M. Junge and Q. Xu, \emph{Noncommutative
Burkholder/Rosenthal inequalities}, Ann. Probab. \textbf{31}
(2003), 948-995.

\bibitem {KP} M.I. Kadec and A. Pe{\l}czy{\'n}ski, \emph{Bases,
lacunary sequences and complemented subspaces in the spaces
${L}\sb{p}$}, Studia Math. \textbf{21} (1961/1962), 161-176.

\bibitem {KR1} R.V. Kadison and J.R. Ringrose, Fundamentals
of the theory of operator algebras II. Academic Press, Orlando,
FL, 1986, Advanced theory.

\bibitem {LT} J. Lindenstrauss and L. Tzafriri, Classical
Banach spaces II. Springer-Verlag, Berlin, 1979, Function spaces.

\bibitem {MT} T. Mart\'{\i}nez and J.L. Torrea, \emph{Operator-valued
martingale transforms}, Tohoku Math. J. \textbf{52} (2000),
449-474.

\bibitem{M} M. Musat, \emph{Interpolation between non-commutative
BMO and non-commutative $L_p$-spaces}, J. Funct. Anal.
\textbf{202} (2003), 195-225.

\bibitem{N} E. Nelson, \emph{Notes on non-commutative integration},
J. Funct. Anal. \textbf{15} (1974), 103--116.

\bibitem {PX1} G. Pisier and Q. Xu, \emph{Non-commutative
martingale inequalities}, Comm. Math. Phys. \textbf{189} (1997),
667-698.

\bibitem {PX2} G. Pisier and Q. Xu, \emph{Non-commutative
$L_p$-spaces}. Handbook of the Geometry of Banach Spaces II (Ed.
W.B. Johnson and J. Lindenstrauss) North-Holland (2003),
1459-1517.

\bibitem {Ran10} N. Randrianantoanina, \emph{Sequences in
non-commutative ${L}\sp{p}$-spaces}, J. Operator Theory
\textbf{48} (2002), 255-272.

\bibitem {R} N. Randrianantoanina, \emph{Non-commutative
martingale transforms}, J. Funct. Anal. \textbf{194} (2002),
181-212.

\bibitem {R3} N. Randrianantoanina, \emph{A weak type inequality
for non-commutative martingales and applications}. Preprint 2004.

\bibitem {T} M. Takesaki, Theory of Operator Algebras I.
Springer-Verlag, New York, 1979.

\bibitem{Tomi} J. Tomiyama, \emph{On the projection of norm one
in ${W}\sp{\ast} $-algebras}, Proc. Japan Acad. \textbf{33}
(1957), 608--612.

\bibitem {X2} Q. Xu, \emph{Analytic functions with values in
lattices and symmetric spaces of measurable operators}, Math.
Proc. Cambridge Philos. Soc. \textbf{109} (1991), 541-563.

\bibitem {X} Q. Xu, \emph{Recent devepolment on non-commutative
martingale inequalities}. Functional Space Theory and its
Applications. Proceedings of International Conference \& 13th
Academic Symposium in China. Ed. Research Information Ltd UK.
Wuhan 2003, 283-314.
\end{thebibliography}

\end{document}